\documentclass{article}
\usepackage[a4paper, margin=2.5cm]{geometry}
\usepackage[title]{appendix}

\usepackage{hyperref}
\usepackage{algorithm}
\usepackage[noend]{algpseudocode}
\usepackage{amsmath}
\usepackage{amssymb}
\usepackage{cleveref}
\usepackage{threeparttable}
\usepackage{booktabs}
\usepackage{subcaption}
\usepackage{tikz}
  \usetikzlibrary{fit}
  \usetikzlibrary{calc}
  \usetikzlibrary{shapes}
  \usetikzlibrary{arrows}
  \usetikzlibrary{positioning}
  \usetikzlibrary{chains}
\usepackage{pgfplots}
\usepackage{pgfplotsthemetol}

\usepackage{natbib}
\usepackage{enumerate}

\newcommand{\F}{\mathcal{F}}
\newcommand{\Ord}{\mathcal{O}}
\newcommand{\curr}{\operatorname{curr}}
\newcommand{\Pk}{\mathcal{P}}
\newcommand{\Dl}{\mathcal{D}}

\usepackage{authblk}

\title{A cost function approximation method for dynamic vehicle routing with docking and LIFO constraints}

\author{Mark\'o Horv\'ath}
\author{Tam\'as Kis}
\author{P\'eter Gy\"orgyi}
\affil{\small HUN-REN Institute for Computer Science and Control, H-1111 Budapest, Kende u. 13-17.}
\affil[ ]{marko.horvath@sztaki.hu, tamas.kis@sztaki.hu, peter.gyorgyi@sztaki.hu}
\begin{document}


\maketitle

\begin{abstract}
In this paper, we study a dynamic pickup and delivery problem with docking constraints.
There is a homogeneous fleet of vehicles to serve pickup-and-delivery requests at given locations.
The vehicles can be loaded up to their capacity, while unloading  has to follow the last-in-first-out (LIFO) rule.
The locations have a limited number of docking ports for loading and unloading, which may force the vehicles to wait.
The problem is dynamic since the transportation requests arrive real-time, over the day. Accordingly, the routes of the vehicles are to be determined dynamically.
The goal is to satisfy all the requests such that a combination of tardiness penalties and traveling costs is minimized.
We propose a cost function approximation based solution method.
In each decision epoch, we solve the respective optimization problem with a perturbed objective function  to ensure the solutions remain adaptable to accommodate new requests. We penalize waiting times and idle vehicles.
We propose a variable neighborhood search based method for solving the optimization problems, and we apply two existing local search operators, and we also introduce a new one.
We evaluate our method using a widely adopted benchmark dataset, and the results demonstrate that our approach significantly surpasses the current state-of-the-art methods.\\

\end{abstract}


\section{Introduction}\label{sec:intro}
Dynamic vehicle routing problems (DVRPs) constitute a rapidly developing field of transportation research, which is certified by a series of recent review papers, e.g., \citet{berbeglia2010dynamic, pillac2013review, bekta2014chapter, psaraftis2016dynamic, soeffker2022stochastic, zhang2022dynamic}.
The growing interest is due to the wide range of real-world application areas such as transportation of goods and people,  services, etc.~\citep{rios2021recent}.

In this paper, we study a DVRP motivated by a real-life problem proposed by Huawei Technologies Co.~Ltd.
A fleet of homogeneous vehicles has to serve pickup-and-delivery requests which occur at given locations over the day. Each request is characterized by a size, a pickup and a delivery location, a release time, and a due date.
The vehicles can be loaded up to their capacity, while unloading  has to follow the last-in-first-out (LIFO) rule.
The locations have a limited number of docking ports for loading and unloading, which may force the vehicles to wait.
The problem is dynamic since the transportation requests arrive real-time, over the day. Accordingly, the routes of the vehicles are to be determined dynamically.
The goal is to satisfy all the requests such that a combination of tardiness penalties and traveling costs is minimized.

An important feature of our problem is that although the pickup and delivery locations are known in advance, the characteristics of the transportation requests are unknown until their release times. Moreover, the distribution of the requests both in time and space may vary over the day. 
Consequently, a single problem instance does not provide exploitable stochastic information. 
However, statistical data may be collected over longer time periods that potentially could be used in solution approaches.

Our problem has two important  constraints that frequently occur  in practice. 
The first one limits the number of the docking ports at each location.
The number of the docking ports is limited by the size of the buildings (warehouses or factories) and also, it may not be economical to have sufficient  crew to serve (load or unload) any number of vehicles simultaneously.
Usually, vehicles do not arrive evenly at a location over the day, which means that in peak periods, some of the vehicles have to wait until a docking port becomes free.  
In several applications the service time of a vehicle may be comparable to the average travel time between two locations, which may lead to high waiting times in case of superficial planning. 
Despite its high practical relevance, the literature on pickup and delivery problems with docking constraints at the locations is rather scarce, we refer to e.g., \citet{cai2022variable, cai2023survey,du2023hierarchical}. 
Most results are for outbound and inbound transportation problems with a depot, see \Cref{sec:constraints} for references.
The second side constraint is the LIFO rule.
This policy means that the last loaded order must be unloaded first.
This rule is necessary, for example, if the vehicles used for transportation have only a single access door for loading and unloading.
Also, if the orders are hazardous, weighty, or fragile, the load rearrangement on the vehicle may consume much time and increase handling costs.
The related literature is summarized in Section~\ref{sec:constraints}.

\paragraph{Main contributions}
In this paper, we propose a cost function approximation  method for the problem at hand, where we manipulate the cost function by introducing two penalty terms.
One of them is directly related to the waiting caused by the docking constraints.
By explicitly penalizing waiting, we expect to make solutions flexible to accommodate new requests in the future.
We model the problem as a sequential decision process, where in each decision epoch a variable neighborhood search (VNS) based method is used to solve the respective optimization problem.
In the VNS method we use two old and a new neighborhood operator.

As a case study, we evaluate our solution approach on a dynamic pickup and delivery problem.
This problem was introduced in \emph{The Dynamic Pickup and Delivery Problem challenge} \citep{hao2022introduction}, hosted by the International Conference on Automated Planning and Scheduling in~2021\footnote{\url{https://icaps21.icaps-conference.org/Competitions/}} (ICAPS 2021).

To sum up, our contributions are as follows.
\begin{itemize}
\item We propose a cost function approximation based method for solving the respective optimization problem in each decision epoch.
One of the penalty terms is directly related to the waiting times caused by the docking constraints.
\item We evaluate our solution procedure on the benchmark instances of the ICAPS 2021 DPDP competition.
The computational experiments show that our method significantly outperforms the state-of-the art methods on this dataset, especially on the large-size instances.
The average improvement on the full dataset is more than 50\% when compared to the best published methods.
\item We demonstrate the benefit of using the suggested penalty terms in separate experiments.
We also explain the mechanism how penalizing waiting improves the solution.
\end{itemize}

According to our best knowledge, no earlier paper applied cost function approximation for a dynamic pickup and delivery  problem with docking constraints and LIFO rule. 

\paragraph{Structure of the paper}
In \Cref{sec:litrev}, we overview the dynamic vehicle routing problems and the related literature.
In \Cref{sec:prob}, we give a formal description of the problem studied along with modeling as a sequential decision process.
In \Cref{sec:sol}, we describe our solution approach and 
in \Cref{sec:comp}, we present our computational results.
Finally, we conclude the paper in \Cref{sec:conc}.


\section{Literature review}\label{sec:litrev}

In this section, we briefly overview the related literature.
In \Cref{sec:pre:vrp}, we narrow our scope to the problem class studied in this paper. 
Modeling with sequential decision process and related solution methods are  summarized in \Cref{sec:pre:seq}.
Finally, \Cref{sec:constraints} is concerned with the LIFO loading rule, and problems with docking constraints.

\subsection{Classification of our vehicle routing problem}\label{sec:pre:vrp}

The problem studied in this paper is \emph{dynamic}, since the input of the problem is received and updated concurrently with the determination of the routes \citep{psaraftis1980dynamic}.
By \citet{pillac2013review}, a dynamic problem is \emph{stochastic}, if some exploitable stochastic knowledge is available on the dynamically revealed information, and \emph{deterministic} otherwise.
In a single instance of our problem, no explicit stochastic information is available, so in this sense our problem is deterministic.

In the \emph{general pickup and delivery problem} a set of routes has to be constructed in order to satisfy transportation request by a fleet of vehicles \citep{savelsbergh1995general}. 
Each transportation request contains a set of pickup locations alongside load quantities and a set of delivery locations along with quantities to unload.
Each request has to be fulfilled by a vehicle without transshipment and exceeding its capacity.
A special case is the \emph{pickup and delivery problem}, where each request has a single pickup location  and a single delivery location. We refer to
\citep{berbeglia2010dynamic} for a survey. 

We conclude that the problem studied in this paper falls in the broad class of DVRPs, and in particular, it is a deterministic dynamic pickup and delivery problem.

\subsection{Modeling and solving DVRPs}\label{sec:pre:seq}
A dynamic vehicle routing problem can be formulated as a \emph{sequential decision process} (SDP). 
An SDP goes through a sequence of \emph{states} that describe the status of the system at distinct time points. 
When passing to the next state, the transition is determined by the dynamic information revealed and the decisions made since the occurrence of the last state. 
The transition may be deterministic or stochastic, in the latter case we have a \emph{Markov decision process}. 
We refer the reader to  \citep{powell2011approximate} for a general introduction to SDPs.

When modeling a~DVRP, decision epochs may occur regularly, or at some special occasions.
The states encode all the information needed to make decisions, and evaluate solution alternatives.
For instance, they may contain the position and current tasks of the vehicles, the set of open service requests, etc.
The decisions are concerned with the assignment of new service requests and possibly  routes to the vehicles.
The transition to the next state determines the new position of the vehicles, or the progress of loading and unloading.
The dynamic information may influence the transition to the new state, e.g., new service requests manifest.

A number of methods have been proposed for solving DVRPs.
The simplest ones are \emph{myopic}, also termed \emph{rolling horizon re-optimization (RO)} methods, that focus only on the current state  without any considerations of future uncertainties. 
Essentially, a series of static problems are solved with the assumption that the realizations will remain unchanged in the future, e.g., no new orders will be requested.
\citet{gendreau1999parallel} solve a vehicle routing problem with time windows with this approach, while \citet{ichoua2000diversion} consider the same problem, but do not allow to change the destination of the moving vehicles.
However, optimizing without acknowledging the future can be counterproductive as it often leads to inflexible solutions.

To alleviate the short-sightedness of the myopic strategies, \citet{mitrovic2004double} propose waiting strategies for a dynamic problem with time windows to delay the dispatch of new orders, and to make decision together with the orders that may be requested in the near future.
\citet{van2004dynamic} consider a dynamic problem of collecting loads, and inserted fictive, anticipated loads into the routes in order to
encourage vehicles to explore fruitful regions (i.e., regions that have a high potential of generating loads).

\emph{Reinforcement learning} for solving DVRPs has gained an increasing attention, we refer to \citep{hildebrandt2023opportunities} for an overview. Below we summarize the most common techniques. 
\emph{Policy function approximation (PFA)} 
is a function that assigns an action to any state,
without any further optimization and without using any forecast of the future information.
\citet{ulmer2019same} apply this method for a problem where parcel pickup stations and autonomous vehicles are combined for same-day delivery.
\citet{ghiani2022scalable}  tackle a pickup and delivery problem using an anticipatory policy by creating priority classes for the requests and their algorithm utilizes a parametric policy function approximation.
\emph{Value function approximation (VFA)} is used to
estimate  the value function, which provides the expected cumulative reward of any given state.
\citet{ulmer2018value} model a multi-period VRP by a Markov decision process, where the set of requests to be accepted in each period is determined by a suitable value function approximation computed offline.
\citet{van2019delivery} consider a delivery dispatching problem and provided a method that  can handle large instances.
For further examples, we refer to the supplementary material of \citet{zhang2022review}.

Beyond reinforcement learning, 
\emph{cost function approximation (CFA)} is a further technique for solving DVRPs.
The key idea is that the problem to be solved in each decision epoch is modified, that is, either the cost function or the problem constraints are slightly perturbed. 
For instance, \citet{riley2019column, riley2020real} consider a dial-a-ride problem with the aim of minimizing the total waiting time.
The authors introduce an extra penalty term for unserved request to ensure that all riders are served in reasonable time.
The penalty associated with a request is increased after each epoch in which the request is not served. \citet{ulmer2020binary} study a retail distribution problem, where familiarity with the delivery location can save service time for the driver, and a related quantity is added to the original objective function.
\citet{hildebrandt2023opportunities} outline an enhancement of CFA methods by first modifying the state received before determining the next action, e.g., by reducing the due dates of some transportation requests, 
or by reducing vehicle capacities, etc.

Other approaches use the stochastic information internally, e.g., via sampling realizations.
For example, \emph{multiple scenario approaches} sample realizations to create a set of scenarios, which are then solved separately, and based on those individual solutions, a consensus decision is made \citep{bent2004scenario, hvattum2006solving, srour2018strategies, dayarian2020crowdshipping}.
For instance, \cite{srour2018strategies} study a stochastic and dynamic pickup and delivery problem with time windows, where the location of future requests is known, but only stochastic information is available about the time windows of future requests.
The authors propose to sample future requests at each decision point and solve a VRP for each scenario, and finally synthesize the next action of the vehicles based on the set of solutions obtained for the different scenarios.
An alternative  approach to sampling scenarios is suggested by \cite{gyorgyi2019probabilistic} for solving the same problem. The stochastic information is used to set up a minimum-cost network flow problem in which source-to-sink routes correspond to vehicle routes and the edges are weighted with  conditional expected values.

\subsection{Docking and LIFO constraints}\label{sec:constraints}
The considered pickup and delivery problem has two important side constraints, namely, the docking constraints, and the LIFO constraints. 
They are of great practical importance and have been studied by the research community.

\paragraph{Docking constraints}
Several papers deal with vehicle routing problems with service restrictions, however, these restrictions mostly apply to a single depot.
In \citet{dabia2019cover} shifts with limited loading capacities are considered.
\citet{hempsch2008vehicle} consider a problem with sorting capacity constraints at the depot and use a local search based method to solve it.
\citet{gromicho2012vehicle} apply a method that uses a  decomposition scheme where columns are generated by a routine based on dynamic programming to solve a variant with limited number of loading docks or limited size of loading crew. 
\citet{van2017vehicle} considers also a vehicle routing problem with a limited number of loading docking ports and propose a method that shifts the start times of the routes in a smart way.

\paragraph{LIFO constraints in pickup and delivery problems}
\citet{cordeau2010branch} propose a branch-and-cut algorithm for a pickup and delivery problem with a single vehicle.
\citet{benavent2015multiple} examine a variant of the problem with multiple vehicles and a special time constraint, and propose a branch-and-cut and a tabu search algorithm as solution methods.
\citet{xu2023dynamic} construct a multi-objective mathematical model for a pickup and delivery problem with transshipments, and propose a method that generates an initial solution by Clarke-Wright saving algorithm and uses both neighborhood search and Q-learning to improve the solution. 
In \citep{carrabs2007variable}, several local search operators (e.g., \emph{couple-exchange}, \emph{block-exchange}, \emph{relocate-couple}, \emph{relocate-block}, \emph{multi-relocate}, \emph{2-opt-L}, \emph{double-bridge}) are proposed to problems with LIFO constraints.

\section{Problem statement}\label{sec:prob}
In the following, we define the problem in detail.
We first define the basic data of the problem, then we present the dynamic problem as a sequential decision process.

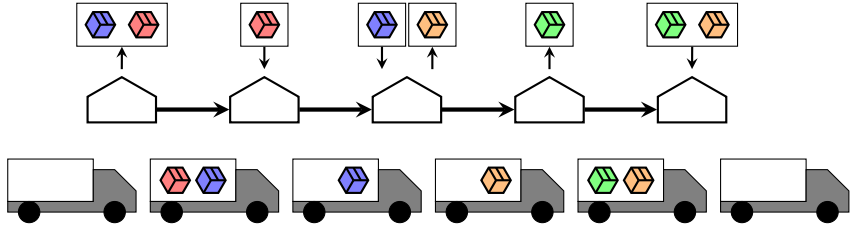
\begin{figure}
\begin{tikzpicture}[scale=0.75,transform shape]
\tikzset{
  dpdpbox/.style n args={2}{
    regular polygon, 
    regular polygon sides=6,
    minimum size=#1,
    rotate=0,
    append after command={
      \pgfextra{
        \draw[thick,fill=#2]
        (\tikzlastnode.corner 1)
        -- (\tikzlastnode.corner 2)
        -- (\tikzlastnode.corner 3)
        -- (\tikzlastnode.corner 4)
        -- (\tikzlastnode.corner 5)
        -- (\tikzlastnode.corner 6)
        -- cycle;
        
        \foreach \i in {2,4,6}{
          \draw[thick] (\tikzlastnode.center) -- (\tikzlastnode.corner \i);
        }
        
        \draw[thick] ($(\tikzlastnode.corner 1)!0.5!(\tikzlastnode.corner 2)$) -- ($(\tikzlastnode.center)!0.5!(\tikzlastnode.corner 6)$);
      }
    }
  },
  dpdpbox/.default={2cm,white}
 }
 
\tikzset{
  dpdpfactory/.style={
    draw,
    thick,
    single arrow, 
    shape border rotate=90, 
    single arrow head extend=0pt, 
    single arrow tip angle=120,
    minimum width=1.2cm,
    minimum height=0.8cm
  }
}

\tikzset{
  pics/dpdpvehicle/.style={
    code={
      \coordinate (zero) at (0,0);
      
      \path[draw,fill=white]
        (zero)++(-0.5*#1,.25*#1)
        -- ++(#1,0)
        -- ++(0,-.5*#1)
        -- ++(-#1,0)
        -- cycle;
        
      \path[draw, fill=black!50!white]
        (zero)++(0.5*#1,0.125*#1)
        -- ++(0.25*#1,0)
        -- ++(0.25*#1,-0.25*#1)
        -- ++(0,-0.25*#1)
        -- ++(-1.5*#1,0)
        -- ++(0,0.125*#1)
        -- ++(1*#1,0)
        -- cycle;
        
      \draw[fill=black] (zero)++(-0.25*#1,-0.375*#1) circle (0.125*#1);
      \draw[fill=black] (zero)++( 0.75*#1,-0.375*#1) circle (0.125*#1);
    }
  }
}

\def\sc{1.25}

\node[dpdpfactory] (f1) at ( 2*\sc,0) {};
\node[dpdpfactory] (f2) at ( 4*\sc,0) {};
\node[dpdpfactory] (f3) at ( 6*\sc,0) {};
\node[dpdpfactory] (f4) at ( 8*\sc,0) {};
\node[dpdpfactory] (f5) at (10*\sc,0) {};

\draw[->,-stealth,ultra thick] (f1) -- (f2);
\draw[->,-stealth,ultra thick] (f2) -- (f3);
\draw[->,-stealth,ultra thick] (f3) -- (f4);
\draw[->,-stealth,ultra thick] (f4) -- (f5);

\begin{scope}[shift={(1.2*\sc,1.5)}]
\node[dpdpbox={0.5cm}{blue!50!white}] (o2p) at (0.5*\sc,0) {};
\node[dpdpbox={0.5cm}{red!50!white}, right=0.5cm] (o1p) at (o2p) {};
\node[draw,fit=(o2p) (o1p)] (p1) {};
\draw[->,thick,-stealth] (f1.north)++(0,0.125) -- (p1);
\end{scope}
\begin{scope}[shift={(4.0*\sc,1.5)}]
\node[dpdpbox={0.5cm}{red!50!white}] (o2d) at (0*\sc,0) {};
\node[draw,fit=(o2d)] (p2) {};
\draw[->,thick,stealth-] (f2.north)++(0,0.125) --(p2);
\end{scope}
\begin{scope}[shift={(5.15*\sc,1.5)}]
\node[dpdpbox={0.5cm}{blue!50!white}] (o4d) at (0.5*\sc,0) {};
\node[dpdpbox={0.5cm}{orange!50!white}, right=0.5cm*\sc] (o3p) at (o4d) {};
\node[draw,fit=(o4d)] (p3d) {};
\node[draw,fit=(o3p)] (p3p) {};
\draw[->,thick,stealth-] (f3.north-|p3d)++(0,0.125) -- (p3d);
\draw[->,thick,-stealth] (f3.north-|p3p)++(0,0.125) -- (p3p);
\end{scope}
\begin{scope}[shift={(8*\sc,1.5)}]
\node[dpdpbox={0.5cm}{green!50!white}] (o4p) at (0*\sc,0) {};
\node[draw,fit=(o4p)] (p4) {};
\draw[->,thick,-stealth] (f4.north)++(0,0.125) -- (p4);
\end{scope}
\begin{scope}[shift={(9.2*\sc,1.5)}]
\node[dpdpbox={0.5cm}{green!50!white}] (o4d) at (0.5*\sc,0) {};
\node[dpdpbox={0.5cm}{orange!50!white}, right=0.5cm] (o3d) at (o4d) {};
\node[draw,fit=(o4d) (o3d)] (p5) {};
\draw[->,thick,stealth-] (f5.north)++(0,0.125) -- (p5);
\end{scope}

\pic at (1*\sc,-1*\sc) {dpdpvehicle=1.5};

\pic at (3*\sc,-1*\sc) {dpdpvehicle=1.5};
\node[dpdpbox={0.5cm}{red!50!white}]  at (2.75*\sc,-1*\sc) {};
\node[dpdpbox={0.5cm}{blue!50!white}] at (3.25*\sc,-1*\sc) {};

\pic at (5*\sc,-1*\sc) {dpdpvehicle=1.5};
\node[dpdpbox={0.5cm}{blue!50!white}] at (5.25*\sc,-1*\sc) {};

\pic at (7*\sc,-1*\sc) {dpdpvehicle=1.5};
\node[dpdpbox={0.5cm}{orange!50!white}] at (7.25*\sc,-1*\sc) {};

\pic at (9*\sc,-1*\sc) {dpdpvehicle=1.5};
\node[dpdpbox={0.5cm}{green!50!white}] at (8.75*\sc,-1*\sc) {};
\node[dpdpbox={0.5cm}{orange!50!white}] at (9.25*\sc,-1*\sc) {};

\pic at (11*\sc,-1*\sc) {dpdpvehicle=1.5};
\end{tikzpicture}
\caption{
An example for loading and unloading during a vehicle route.
Orders are depicted as boxes, factories are depicted as pentagons.
Unloading the orders from a vehicle has to follow the LIFO rule, see the order of orders above the factories, and also the position of the orders on the vehicle.
}
\label{fig:loading}
\end{figure}

Given a finite set of factories $\F$, a finite set of orders $\Ord$, and a fleet of homogeneous vehicles $\mathcal{V}$.
The \emph{distance} and the \emph{travel time} between factories~$f_i, f_j \in \F$ are denoted with~$\operatorname{dist}(f_i,f_j)$ and~$\operatorname{travel}(f_i,f_j)$, respectively.
Each order~$o_i\in\Ord$ is described by a tuple $(f^p_i,f^d_i,t^p_i,t^d_i,q_i,h^p_i,h^d_i)$, where $f^p_i$ and $f^d_i$ represent the \emph{pickup factory} and the \emph{delivery factory}, respectively, $t^p_i$ is the \emph{release time}, $t^d_i$ is the \emph{due date}, $q_i$ is the order \emph{quantity}, and $h^p_i$ and $h^d_i$ are the times required to \emph{load} and to \emph{unload} the order, respectively. 
Order~$o_i$ becomes known only at the release time~$t^p_i$.
The total quantity of those orders that can be carried by a vehicle at any given moment is  limited by a constant~$Q$.
Initially, each vehicle is empty and parks at a given factory.
Unloading the orders from a vehicle has to follow the LIFO rule, refer to \Cref{fig:loading} for an example.

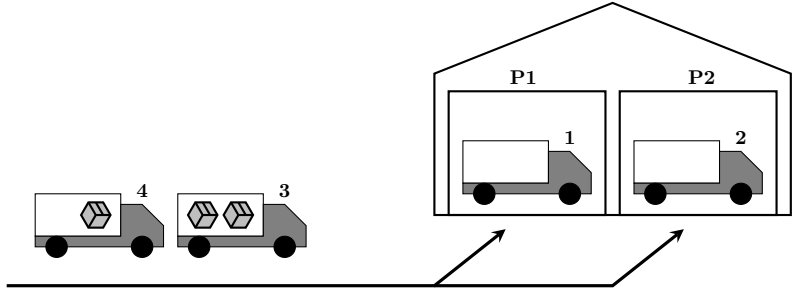
\begin{figure}
\begin{tikzpicture}[scale=0.75,transform shape]
\tikzset{
  dpdpbox/.style n args={2}{
    regular polygon, 
    regular polygon sides=6,
    minimum size=#1,
    rotate=0,
    append after command={
      \pgfextra{
        \draw[thick,fill=#2]
        (\tikzlastnode.corner 1)
        -- (\tikzlastnode.corner 2)
        -- (\tikzlastnode.corner 3)
        -- (\tikzlastnode.corner 4)
        -- (\tikzlastnode.corner 5)
        -- (\tikzlastnode.corner 6)
        -- cycle;
        
        \foreach \i in {2,4,6}{
          \draw[thick] (\tikzlastnode.center) -- (\tikzlastnode.corner \i);
        }
        
        \draw[thick] ($(\tikzlastnode.corner 1)!0.5!(\tikzlastnode.corner 2)$) -- ($(\tikzlastnode.center)!0.5!(\tikzlastnode.corner 6)$);
      }
    }
  },
  dpdpbox/.default={2cm,white}
 }

\tikzset{
  pics/dpdpvehicle/.style 2 args={
    code={
      \coordinate (zero) at (0,0);
      
      \path[draw,fill=white]
        (zero)++(-0.5*#1,.25*#1)
        -- ++(#1,0)
        -- ++(0,-.5*#1)
        -- ++(-#1,0)
        -- cycle;
        
      \path[draw, fill=black!50!white]
        (zero)++(0.5*#1,0.125*#1)
        -- ++(0.25*#1,0)
        -- ++(0.25*#1,-0.25*#1)
        -- ++(0,-0.25*#1)
        -- ++(-1.5*#1,0)
        -- ++(0,0.125*#1)
        -- ++(1*#1,0)
        -- cycle;
        
      \draw[fill=black] (zero)++(-0.25*#1,-0.375*#1) circle (0.125*#1);
      \draw[fill=black] (zero)++( 0.75*#1,-0.375*#1) circle (0.125*#1);
      
      \node[above] at (0.75*#1,0.125*#1) {\textbf{#2}};
    }
  }
}

\def\sc{1.25}

\begin{scope}[shift={(5*\sc,0)}]
\path[draw,thick] (0*\sc,0*\sc) -- (0*\sc,2*\sc) -- (2.5*\sc,3*\sc) -- (5*\sc,2*\sc) -- (5*\sc,0*\sc) -- cycle;
\draw[thick] (0.2*\sc,0*\sc) rectangle (2.4*\sc,1.75*\sc);
\draw[thick] (2.6*\sc,0*\sc) rectangle (4.8*\sc,1.75*\sc);
\pic at (1*\sc,0.75*\sc) {dpdpvehicle={1.5}{1}};
\pic at (3.4*\sc,0.75*\sc) {dpdpvehicle={1.5}{2}};

\node[above] at (1.25*\sc,1.75*\sc) {\textbf{P1}};
\node[above] at (3.75*\sc,1.75*\sc) {\textbf{P2}};
\end{scope}

\pic at (0*\sc,0*\sc) {dpdpvehicle={1.5}{4}};
\node[dpdpbox={0.5cm}{black!25!white}]  at (0.25*\sc,0*\sc) {};

\pic at (2*\sc,0*\sc) {dpdpvehicle={1.5}{3}};
\node[dpdpbox={0.5cm}{black!25!white}]  at (1.75*\sc,0*\sc) {};
\node[dpdpbox={0.5cm}{black!25!white}] at (2.25*\sc,0*\sc) {};

\path[draw,very thick,->,-stealth] (-1*\sc,-1*\sc) -- (5*\sc,-1*\sc) -- (6*\sc,-0.2*\sc);
\path[draw,very thick,->,-stealth]                    (5*\sc,-1*\sc) -- (7.5*\sc,-1*\sc) -- (8.5*\sc,-0.2*\sc);
\end{tikzpicture}
\caption{
An example for serving vehicles at a factory.
The two docking ports, P1 and P2, of the factory are occupied by vehicles~1 and~2, respectively, thus the service of vehicles~3 and~4 cannot start at this moment, and the vehicles must wait until a docking port becomes free.
}
\label{fig:docking}
\end{figure}

Each factory has a given number of docking ports for loading and unloading.
Vehicles are served on a first-come-first-served basis, that is, if a vehicle arrives at a factory and all ports are occupied, the service of the vehicle cannot begin immediately, but the vehicle has to wait until one of the docking ports becomes free, and no vehicle that arrived earlier is waiting for a port. This is illustrated in \Cref{fig:docking}. 
The time elapsed between the arrival of the vehicle and the start of service is called the \emph{waiting time}.
Serving a vehicle decomposes to dock approaching, then unloading some carried orders, and finally loading some new orders.
The \emph{service time} of a vehicle is the sum of the factory-independent \emph{dock approaching time}, $h^{\operatorname{docking}}$, and the sum of the unloading and loading times of the corresponding orders.
For more details we refer to Appendix~\ref{sec:apx:serving}.
After serving a vehicle, the port becomes free, and the vehicle may park at the factory, or depart to the next factory on its route.

The goal is to route the vehicles so that all orders are served, and the weighted sum of the total distance traveled and the total tardiness of the orders is minimized.

\subsection{Modeling as a sequential decision process}\label{sec:SDP}
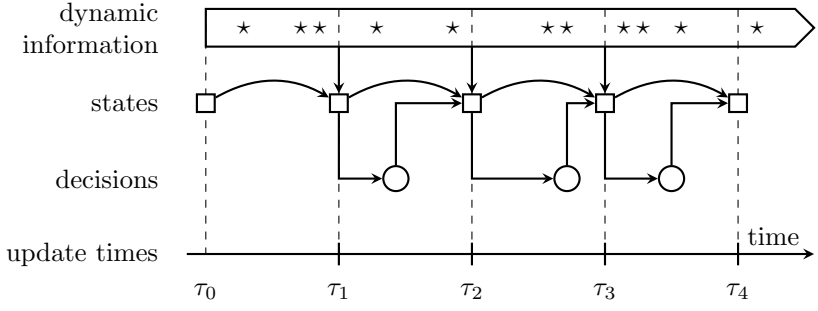
\begin{figure}
\begin{tikzpicture}
\tikzstyle{nrealization} = [thick,shape=signal,signal from=west, signal to=east,
     align=center,draw=black,font=\sffamily,on chain,minimum height=0.5cm,
     inner xsep=1em]
     
\tikzstyle{nstate} = [draw, thick, rectangle, fill= white, minimum height= 0.15cm, minimum width= 0.15cm]
\tikzstyle{ndecision} = [draw, thick, circle, fill= white, minimum height=0.1cm]
\tikzstyle{connection} = [->, -stealth, thick]

\draw[dashed] (0.00,0.25) -- (0.00,-3);
\draw[dashed] (1.75,0.25) -- (1.75,-3);
\draw[dashed] (3.50,0.25) -- (3.50,-3);
\draw[dashed] (5.25,0.25) -- (5.25,-3);
\draw[dashed] (7.00,0.25) -- (7.00,-3);
\begin{scope}[]
 \path[start chain=going right,node distance=0ex]
    node[nrealization,signal from=nowhere,minimum width= 8.00 cm] at (3.875,0) {};
    
\node at (0.50,0) {$\star$};
\node at (1.25,0) {$\star$};
\node at (1.50,0) {$\star$};
\node at (2.25,0) {$\star$};
\node at (3.25,0) {$\star$};
\node at (4.50,0) {$\star$};
\node at (4.75,0) {$\star$};
\node at (5.50,0) {$\star$};
\node at (5.75,0) {$\star$};
\node at (6.25,0) {$\star$};
\node at (7.25,0) {$\star$};
\end{scope}
\begin{scope}[shift={(0,-1)}]
\node[nstate] (s0) at (0.00,0) {};
\node[nstate] (s1) at (1.75,0) {};
\node[nstate] (s2) at (3.50,0) {};
\node[nstate] (s3) at (5.25,0) {};
\node[nstate] (s4) at (7.00,0) {};
\draw[connection] (s0) to [out=30,in=150] (s1);
\draw[connection] (s1) to [out=30,in=150] (s2);
\draw[connection] (s2) to [out=30,in=150] (s3);
\draw[connection] (s3) to [out=30,in=150] (s4);
\end{scope}
\begin{scope}[shift={(0,-2)}]
\node[ndecision] (d1) at (2.50,0) {};
\node[ndecision] (d2) at (4.75,0) {};
\node[ndecision] (d3) at (6.125,0) {};
\end{scope}
\begin{scope}[shift={(0,-3)}]
\draw[thick,->,-stealth] (-0.25,0) -- (8.00,0);
\draw[thick] (1.75,-0.15) -- (1.75,0.15);
\draw[thick] (3.50,-0.15) -- (3.50,0.15);
\draw[thick] (5.25,-0.15) -- (5.25,0.15);
\draw[thick] (7.00,-0.15) -- (7.00,0.15);
\node at (0,-0.5) {$\tau_0$};
\node at (1.75,-0.5) {$\tau_1$};
\node at (3.50,-0.5) {$\tau_2$};
\node at (5.25,-0.5) {$\tau_3$};
\node at (7.00,-0.5) {$\tau_4$};
\node[above right] at (7,0) {time};
\end{scope}
\draw[->, -stealth, thick, to path={|- (\tikztotarget)}]
  (s1) edge (d1)
  (d1) edge (s2)
  (s2) edge (d2)
  (d2) edge (s3)
  (s3) edge (d3)
  (d3) edge (s4);
  
\draw[connection] (1.75,-0.25) -- (s1);
\draw[connection] (3.50,-0.25) -- (s2);
\draw[connection] (5.25,-0.25) -- (s3);

\begin{scope}[shift={(-0.5,0)}]
\node[left, align= right] at (0, 0) {dynamic\\information};
\node[left] at (0,-1) {states};
\node[left] at (0,-2) {decisions};
\node[left] at (0,-3) {update times};
\end{scope}
\end{tikzpicture}
\caption{
Sequential decision process.
}
\label{fig:sdp}
\end{figure}

We model our problem as an SDP \citep{powell2011approximate,soeffker2022stochastic}.
The process goes through a sequence of states $s_0, s_1,\ldots$, each corresponding to an instance of the problem's decision model at time points $\tau_0, \tau_1,\ldots$, as illustrated in \Cref{fig:sdp}.
Between two consecutive time points $\tau_k$ and $\tau_{k+1}$, 
an algorithm computes some actions based on the  state  $s_k$  at $\tau_k$, which comprises the new dynamic information revealed between $\tau_{k-1}$ and $\tau_k$.
The actions along with the dynamic information between  $\tau_k$ and $\tau_{k+1}$ determine the transition of the system to the next state $s_{k+1}$ at $\tau_{k+1}$.

\subsubsection{Decision epochs}

The \emph{update times}~$\tau_k$ ($k=0,1,\ldots,K$) divide the operating horizon into \emph{epochs} of length~$\Delta$ each, that is, $\tau_0 = 0$ and $\tau_{k+1} = \tau_k + \Delta$.
The orders are requested in a given finite planning horizon (e.g., one day), however, the operating horizon may be longer if we are not able to complete the orders within the planning horizon.
The decision process ends if all the orders are delivered.

\subsubsection{Dynamic information}
The \emph{realization}~$\omega_k$ of the dynamic information at  update time $\tau_k$ comprises the orders requested in the previous epoch, i.e., $\omega_k = \{ o_i \in \Ord : \tau_{k-1} < t^p_i \leq \tau_k \}$.

\subsubsection{States}

A \emph{state}~$s_k$ at update time $\tau_k$ is a tuple $(\tau_k,\Phi_k,\tilde{\omega}_k)$, where $\Phi_k = \{ \Phi_{k,v} : v\in \mathcal{V} \}$ is the status of the vehicles, and $\tilde{\omega}_k$ is the set of unprocessed orders.

The set $\tilde{\omega}_k$ comprises the dynamic information $\omega_k$, and those  orders which are released not later than $\tau_{k-1}$, but not picked up until $\tau_k$.

At any time moment a vehicle is either located at a factory, which is then its {\em current factory\/} or is on the way to its {\em destination factory\/}.
The \emph{status of a vehicle~$v$\/} is described by a tuple~$\Phi_{k,v} = (\phi^{\curr}_{k,v},\mathcal{C}_{k,v},\theta_{k,v})$, where
\begin{itemize}
\item $\phi^{\curr}_{k,v} = (f^{\curr}_{k,v}, td^{\curr}_{k,v})$  provides the current factory and the  \emph{earliest departure time} if the vehicle is located at the factory $f^{curr}_{k,v}$ at time $\tau_k$, while $\phi^{\curr}_{k,v} = \emptyset$ if the vehicle is on the way to a factory at time $\tau_k$.

\item $\mathcal{C}_{k,v}$ is the list of orders carried by the vehicle, sorted in the order of loading. If $\phi^{\curr}_{k,v} \neq \emptyset$, then 
$\mathcal{C}_{k,v}$ contains all the orders that the vehicle had to pickup at $f^{curr}_{k,v}$, and does not contain those orders that were delivered to that factory.
\item $\theta_{k,v}$ is the route plan of the vehicle consisting of a sequence of tuples each corresponding to a factory:
\[
  \theta_{k,v} = \left((f^j_{k,v}, ta^j_{k,v}, td^j_{k,v}, \mathcal{D}^j_{k,v}, \mathcal{P}^j_{k,v}) : 1\leq j\leq \ell_{k,v}\right).
\]
The $j$-th tuple in this list corresponds to 
 the $j$th factory~$f^j_{k,v}$ of the route with arrival time $ta^j_{k,v}$,  departure time $td^j_{k,v}$, and with the list of orders $\mathcal{D}^j_{k,v}$ to be unloded, and $\mathcal{P}^j_{k,v}$ to be loaded, respectively.
The first factory visited in  $\theta_{k,v}$ is  the \emph{destination factory} of the vehicle.
Each route plan must be  \emph{feasible}, i.e., it has to  fulfill the \emph{fundamental routing constraints}, the \emph{capacity constraints}, and the \emph{LIFO constraints}, refer to Appendix~\ref{sec:apx:feasibility} for full details.
\end{itemize}

In the initial state~$s_0$ at~$\tau_0$, vehicles have only current factories, which coincide with their initial factories.

\subsubsection{Actions}
Actions are taken at update points.
An action is merely a feasible route plan for each vehicle.
That is, at update point $\tau_k$, let $\Theta(s_k)$ be the set of all possible tuples of feasible and mutually compatible route plans for the vehicles from which exactly one tuple $x_k = ( \theta^x_{k,v} : v \in \mathcal{V})$ must be chosen.
Notice that the $\theta^x_{k,v}$ must be feasible route plans, and they must be \emph{mutually compatible}, i.e., the same order cannot be served by route plans of distinct vehicles.

A further constraint is that if $\theta_{k,v}$ is non-empty, then the first factory  visited in $\theta_{k,v}$ and $\theta^x_{k,v}$ must be the same, i.e., $f^1_{k,v} = f^{x,1}_{k,v}$. Moreover, $ta^1_{k,v} = ta^{x,1}_{k,v}$, and $\Dl^1_{k,v} = \Dl^{x,1}_{k,v}$ must hold.
However, the set of orders to pickup and thus the departure time from $f^1_{k,v}$ may be different in $\theta_{k,v}$ and $\theta^x_{k,v}$.

\subsubsection{Reward function}
The reward function assigns a value to a (state, action) pair $(s_k, x_k)$.
Let $\mathbf{f}_1$ denote the total distance traveled by the vehicles, that is,
\[
  \mathbf{f}_1(s_k,x_k) = 
    \sum_{v\in\mathcal{V}}                          \operatorname{dist}(f^{\curr}_{k,v},f^{x,1}_{k,v}) +
    \sum_{v\in\mathcal{V}}\sum_{j=2}^{\ell^x_{k,v}} \operatorname{dist}(f^{x,j-1}_{k,v},f^{x,j}_{k,v}),
\]
where $\operatorname{dist}(f^{\curr}_{k,v},f^{x,1}_{k,v}) = 0$ if $\phi^{\curr}_{k,v} = \emptyset$ or $\ell^x_{k,v} = 0$, 
and let $\mathbf{f}_2$ be the total tardiness anticipated:
\[
  \mathbf{f}_2(s_k,x_k) = \sum_{v\in\mathcal{V}}
  \sum_{j=1}^{\ell^x_{k,v}}\sum_{o_i \in \mathcal{D}^{x,j}_{k,v}} \max(0,ta^{x,j}_{k,v} - t^d_i).
\]
Then, the reward function is
\begin{equation}\label{eq:reward_0}
R_0(s_k, x_k) = \lambda_1 \mathbf{f}_1(s_k,x_k) + \lambda_2 \mathbf{f}_2(s_k,x_k),
\end{equation}
where $\lambda_1, \lambda_2>0$ are appropriate multipliers.
We will modify this function in our CFA method in \Cref{sec:CFA_cost}.

\subsubsection{Transition}\label{sec:prob:transition}
After the selection of action $x_k$, the decision process transitions to the next state~$s_{k+1}$ at update point $\tau_{k+1}$.
The positions of the vehicles, and the lists of carrying orders are updated.
Briefly stated, if a vehicle $v$ arrives at its destination factory before $\tau_{k+1}$, that factory becomes its current factory, and the next factory to visit, if any, becomes the new destination factory. As a result, $\theta^x_{k,v}$ is transformed to the route plan $\theta_{k+1,v}$ in the new state  $s_{k+1}$.
$\mathcal{C}_{k+1,v}$ is obtained from $\mathcal{C}_{k,v}$ by removing the orders delivered at the current factory of vehicle $v$ in state $s_k$, and adding to it the pickup orders, if any.
If vehicle $v$ departed from its current factory between $\tau_k$ and $\tau_{k+1}$, then no current factory will be associated with the vehicle, and the factory to which the vehicle is heading will be the destination factory.
Then $\Phi_{k+1,v} = ( \mathcal{C}_{k+1,v}, \phi^{\curr}_{k+1,v}, \theta_{k+1,v})$ for each $v \in \mathcal{V}$.
For details, we refer to Appendix~\ref{sec:apx:transition}.

\subsubsection{Objective function}\label{sec:prob:dynamic:obj}

The SDP eventually creates a feasible route $\theta_{v}$ for each vehicle $v$, and the routs are mutually compatible and serve all the requests. 
Moreover, we assume that for each $v \in \mathcal{V}$,  $f^1_{v}$ coincides with the initial factory of vehicle $v$ in the initial state $s_0$.
A formal definition of the feasibility of a solution is given in Appendix~\ref{sec:apx:feas_solution}.
Then $x = (\theta_{v}\ :\ v \in \mathcal{V})$ is a solution of the problem.
The solution is evaluated by the cost function
\begin{equation}\label{eq:obj}
  \mathit{cost}(x) = \lambda_1\mathbf{f}_1(s_0,x) + \lambda_2\mathbf{f}_2(s_0,x).
\end{equation}

\subsection{Example}
\begin{figure}
\begin{tikzpicture}
\def\sc{1.5}
\tikzstyle{nfactory} = [fill= black,circle,minimum height= 0.15cm]
\tikzstyle{nhistfactory} = [fill= black!50!,circle,minimum height= 0.15cm]
\tikzstyle{nvehicle} = [draw= black,fill= white,rectangle,minimum width= 0.5cm, minimum height= 0.25cm]
\tikzstyle{nlabel} = [draw= black,fill= white,rectangle,scale=0.8]
\tikzstyle{dhistory} = [draw= black!50!white,thick,->,-stealth]
\tikzstyle{dcurrent} = [draw= black,very thick,->,-stealth]
\tikzstyle{dfuture} = [draw= black,thick,dashed,->,-stealth]
\tikzstyle{nplan} = [fill= white, rectangle,scale=0.5]


\begin{scope}[scale=\sc]
\node[nhistfactory] (f10) at (-0.7,1) {};
\node[above left] at (f10) {$f_0$};
\node[nhistfactory] (f11) at (0.7,0.8) {};
\node[above right] at (f11) {$f_1$};
\node[nhistfactory] (f12) at (1.0,0.2) {};
\node[above right] at (f12) {$f_2$};
\node[nfactory] (f13) at (0,0) {};
\node[above left] at (f13) {$f_3$};
\node[nfactory] (f14) at (-.3,-0.8) {};
\node[below left] at (f14) {$f_4$};
\node[nfactory] (f15) at (0.7,-1) {};
\node[above right] at (f15) {$f_5$};
\node[nfactory,fill= white,draw = black] (f16) at (-1,-0.2) {};
\node[above left] at (f16) {$f_6$};

\draw[dhistory] (f10) -- (f11);
\draw[dhistory] (f11) -- (f12);
\draw[dcurrent] (f12) -- (f13);
\draw[dfuture] (f13) -- (f14);
\draw[dfuture] (f14) -- (f15);

\node[nplan,below= 0.25cm] at       (f10) {$[0,10]$};
\node[nplan,below left= 0.25cm] at  (f11) {$[21,25],o_2^+,o_1^+$};
\node[nplan,below= 0.25cm] at       (f12) {$[36,39],o_1^-$};
\node[nplan,below right= 0.25cm] at (f13) {$[50,53],o_2^-$};
\node[nplan,above right= 0.25cm] at (f14) {$[64,67],o_3^+$};
\node[nplan,below= 0.25cm] at       (f15) {$[78,81],o_3^-$};
\node[nplan,below= 0.25cm] at       (f16) {$o_4\to f_5$};

\node[nvehicle] at ($(f12)!0.5!(f13)$) {};

\node[nlabel,above right] at (-0.9*\sc,-\sc) {$s_k$ ($\tau_k=40$)};
\end{scope}


\draw[] (1.5*\sc,-1.5*\sc) -- (1.5*\sc,1.5*\sc);


\begin{scope}[scale=\sc,shift={(2*\sc,0)}]
\node[nhistfactory] (f10) at (-0.7,1) {};
\node[above left] at (f10) {$f_0$};
\node[nhistfactory] (f11) at (0.7,0.8) {};
\node[above right] at (f11) {$f_1$};
\node[nhistfactory] (f12) at (1.0,0.2) {};
\node[above right] at (f12) {$f_2$};
\node[nfactory] (f13) at (0,0) {};
\node[above left] at (f13) {$f_3$};
\node[nfactory] (f14) at (-.3,-0.8) {};
\node[below left] at (f14) {$f_4$};
\node[nfactory] (f15) at (0.7,-1) {};
\node[above right] at (f15) {$f_5$};
\node[nfactory] (f16) at (-1,-0.2) {};
\node[above left] at (f16) {$f_6$};

\draw[dhistory] (f10) -- (f11);
\draw[dhistory] (f11) -- (f12);
\draw[dcurrent] (f12) -- (f13);
\draw[dfuture] (f13) -- (f16);
\draw[dfuture] (f16) -- (f14);
\draw[dfuture] (f14) -- (f15);

\node[nplan,below= 0.25cm] at       (f10) {$[0,10]$};
\node[nplan,below left= 0.25cm] at  (f11) {$[21,25],o_2^+,o_1^+$};
\node[nplan,below= 0.25cm] at       (f12) {$[36,39],o_1^-$};
\node[nplan,below= 0.25cm] at       (f13) {$[50,53],o_2^-$};
\node[nplan,below= 0.25cm] at (f16) {$[64,67],o_4^+$};
\node[nplan,above right= 0.25cm] at (f14) {$[78,81],o_3^+$};
\node[nplan,below= 0.25cm] at       (f15) {$[92,96],o_3^-,o_4^-$};

\node[nvehicle] at ($(f12)!0.5!(f13)$) {};

\node[nlabel,above right] at (-0.9*\sc,-\sc) {$x_k$};
\end{scope}


\draw[] (4.5*\sc,-1.5*\sc) -- (4.5*\sc,1.5*\sc);


\begin{scope}[scale=\sc,shift={(4*\sc,0)}]
\node[nhistfactory] (f10) at (-0.7,1) {};
\node[above left] at (f10) {$f_0$};

\node[nhistfactory] (f11) at (0.7,0.8) {};
\node[above right] at (f11) {$f_1$};

\node[nhistfactory] (f12) at (1.0,0.2) {};
\node[above right] at (f12) {$f_2$};

\node[nfactory] (f13) at (0,0) {};
\node[above left] at (f13) {$f_3$};

\node[nfactory] (f14) at (-.3,-0.8) {};
\node[below left] at (f14) {$f_4$};

\node[nfactory] (f15) at (0.7,-1) {};
\node[above right] at (f15) {$f_5$};

\node[nfactory] (f16) at (-1,-0.2) {};
\node[above left] at (f16) {$f_6$};

\draw[dhistory] (f10) -- (f11);
\draw[dhistory] (f11) -- (f12);
\draw[dhistory] (f12) -- (f13);
\draw[dcurrent] (f13) -- (f16);
\draw[dfuture] (f16) -- (f14);
\draw[dfuture] (f14) -- (f15);

\node[nplan,below= 0.25cm] at       (f10) {$[0,10]$};
\node[nplan,below left= 0.25cm] at  (f11) {$[21,25],o_2^+,o_1^+$};
\node[nplan,below= 0.25cm] at       (f12) {$[36,39],o_1^-$};
\node[nplan,below= 0.25cm] at       (f13) {$[50,53],o_2^-$};
\node[nplan,below= 0.25cm] at (f16) {$[64,67],o_4^+$};
\node[nplan,above right= 0.25cm] at (f14) {$[78,81],o_3^+$};
\node[nplan,below= 0.25cm] at       (f15) {$[92,96],o_3^-,o_4^-$};

\node[nvehicle] at (f13) {};

\node[nlabel,above right] at (-0.9*\sc,-\sc) {$s_{k+1}$ ($\tau_{k+1} = 50$)};
\end{scope}

\end{tikzpicture}
\caption{
An example for a route plan of a vehicle at different states and actions.
Left pane shows the route at an intermediate state of the decision process.
Center pane shows the updated route according to an action.
Right pane shows the route at the next state of the decision process.
}
\label{fig:ex:states}
\end{figure}
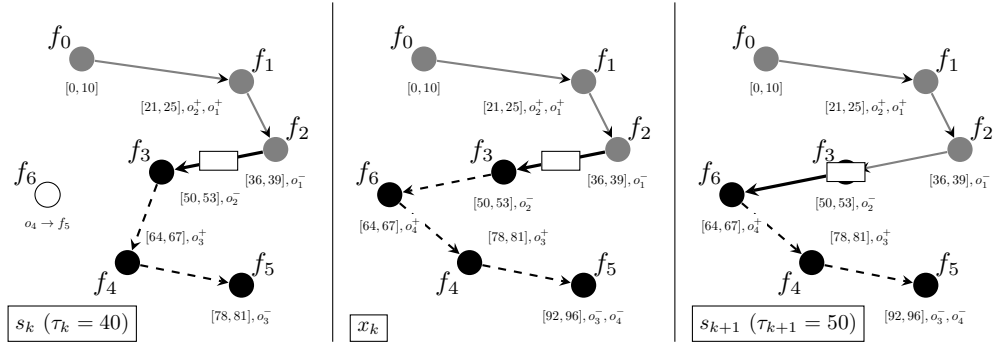

In \Cref{fig:ex:states}, we depict the route of a single vehicle~$v$ at different states.
The gray nodes and edges represent the route of the vehicle before the respective update time points. The black thick edges outline the  route to the destination factory, which cannot be changed, while the route indicated by dashed edges can be modified.
To ease the calculations, we assume that the travel time between any two factories is 11~minutes, the dock approaching time is 2~minutes, and the loading/unloading time for an order is 1~minute.
The length of the epochs is $\Delta =10$ minutes.

Before $\tau_4 = 40$, three orders have arrived: order~$o_1$ from $f_1$ to $f_2$, order~$o_2$ from $f_1$ to $f_3$, and order~$o_3$ from $f_4$ to $f_5$.
The vehicle departed from its initial factory~$f_0$ at time~$10$ and traveled to factory~$f_1$.
The vehicle arrived at factory~$f_1$ at time~$21$, picked up orders $o_2$ and $o_1$ and departed at time~$25$.
The vehicle arrived at factory~$f_2$ at time~$36$, delivered order~$o_1$ and left the factory at time~$39$ to travel to factory~$f_3$.

\paragraph{State~$s_k$ ($\tau_k=40$)}
The left pane of \Cref{fig:ex:states} shows the state~$s_k$ at update point $\tau_k=40$.
The vehicle is currently on the way to factory~$f_3$, where it will arrive at time~$50$.
According to the tentative route plan, the vehicle after that will travel to factory~$f_4$ to pickup and deliver order~$o_3$ to factory~$f_5$.
Thus, the status of the vehicle is given by $\phi^{\curr}_{k,v} = \emptyset$,  $\mathcal{C}_{k,v} = (o_2)$, and $\theta_{k,v} = (\left(f_3,50,53,(o_2),\emptyset)\right.$, $(f_4,64,67,\emptyset,(o_3))$, $\left.(f_5,78,81,(o_3),\emptyset)\right)$.
Order~$o_4$ from factory~$f_6$ to~$f_5$ is also revealed in the previous epoch (see the white circle node), thus $\omega_k = \{o_4\}$.
Since order~$o_1$ is already delivered, and order~$o_2$ is already picked up, $\tilde{\omega}_k = \{o_3,o_4\}$.

\paragraph{Action~$x_k$}
The decision maker decided to insert factory~$f_6$ into the the tentative route plan of the vehicle, that is, $x_k = (\theta^x_{k,v})$, where $\theta^x_{k,v} = \left((f_3,50,53,(o_2),\emptyset)\right.$, $(f_6,64,67,\emptyset,(o_4))$, $(f_4,78,81,\emptyset,(o_3))$, $\left.(f_5,92,96,(o_3,o_4),\emptyset)\right)$.
In the center pane, we depict the updated route plan of the vehicle.

\paragraph{State~$s_{k+1}$ ($\tau_{k+1}=50$)}
The right pane of \Cref{fig:ex:states} shows the state~$s_{k+1}$ at update point $\tau_{k+1}=50$.
Since the vehicle just arrived at factory~$f_3$, it became its current factory.
According to action~$x_k$, after this visit the vehicle will depart toward its new destination factory~$f_6$ to pickup order~$o_4$.
Thus, the status is given by $\mathcal{O}_{k+1,v} = \emptyset$,
$\phi^{\curr}_{k+1,v} = (f_3,53)$,
and $\theta_{k+1,v} = \left((f_6,64,67,\emptyset,(o_4))\right.$, $(f_4,78,81,\emptyset,(o_3))$, $\left.(f_5,92,96,(o_3,o_4),\emptyset)\right)$.
No new orders are requested in the previous epoch, thus~$\omega_{k+1} = \emptyset$, however, the pickup of orders~$o_3$ and $o_4$ can be still changed, thus $\tilde{\omega}_{k+1} = \{o_3,o_4\}$.

\section{CFA approach for solving the routing problem in each epoch}\label{sec:sol}
In this section, we propose a cost function approximation based approach to solve the routing problem in each epoch.
First, we add penalty terms to (\ref{eq:reward_0}) in \Cref{sec:CFA_cost}, then present our VNS procedure in \Cref{sec:method} after describing a new representation of the vehicle routes.

\subsection{Cost function approximation}\label{sec:CFA_cost}
In this section we modify the reward function (\ref{eq:reward_0}) by adding two penalty terms to it. On the one hand, we will penalize waiting for service at the factories, and on the other hand, the idle vehicles.

\paragraph{Penalizing waiting for service}
If the number of vehicles is much larger than the docking capacity of the factories, the vehicles may have to spend a considerable time with queuing. However, waiting times may create large delays in delivery, therefore, it is better to avoid them.
Let $\eta_{k,v}^{x,j}$ be the waiting time (i.e., the time between the arrival and the start of service) at the $j$th  factory visited in the route plan $\theta^x_{k,v}$ of vehicle~$v$.
Then the total waiting time is
\[
  \mathbf{f}_3(s_k,x_k) = \sum_{v\in\mathcal{V}}\sum_{j=1}^{\ell^x_{k,v}} \eta^{x,j}_{k,v}.
\]

\paragraph{Idle vehicles}
We noticed that in some cases the assignment of the first orders significantly affects the subsequent delivery times.
That is, intuitively good solutions (orders with a common pickup factory were assigned to the same vehicle) in the first epochs caused often irreversible tardiness for future orders.
In those cases, it proved better to spread the initial orders between several vehicles, in order to keep as many vehicles moving as possible.
We call a vehicle \emph{idle}, if it has no destination factory, and it will be available in the next epoch.
Then, the total number of idle vehicles is
\[
  \mathbf{f}_4(s_k,x_k) = | \{ v\in\mathcal{V} : \theta^x_{k,v} = \emptyset\ \text{and}\ td^{\curr}_{k,v} < \tau_{k+1} \} |.
\]
\paragraph{Perturbed reward function} We will use the following reward function
\begin{equation}\label{eq:obj:manipulated}
  R(s_k, x_k) = \lambda_1\mathbf{f}_1(s_k,x_k) + \lambda_2 \mathbf{f}_2(s_k,x_k) + \lambda_3\mathbf{f}_3(s_k,x_k) + \lambda_4\mathbf{f}_4(s_k,x_k),
\end{equation}
where $\lambda_1, \lambda_2, \lambda_3, \lambda_4 \geq 0$ are appropriate multipliers.
The  effect of penalty terms on the efficiency of our method will be investigated in Section~\ref{sec:comp:eval}.

\subsection{New representation of routes and the VNS based procedure}\label{sec:method}
Before delving into the details of our solution procedure, first we introduce the internal representation of vehicle routes.

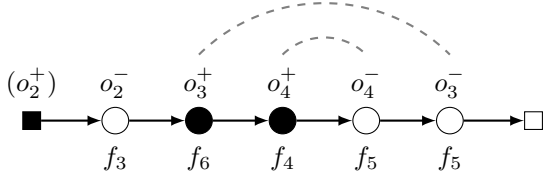
\begin{figure}
\begin{tikzpicture}
\tikzstyle{hnode}=[draw,fill=black,rectangle,minimum height=0.5em, minimum width=0.5em];
\tikzstyle{tnode}=[draw,fill=white,rectangle,minimum height=0.5em, minimum width=0.5em];
\tikzstyle{pnode}=[draw,fill=black,circle,minimum height=1em, minimum width=1em];
\tikzstyle{dnode}=[draw,fill=white,circle,minimum height=1em, minimum width=1em];
\tikzstyle{connection} = [->,thick,latex-,-latex];
\tikzstyle{couple} = [thick,dashed,black!50!white];

\def\spc{1.1}

\foreach \x in {0,...,11}
{
  \coordinate (c\x) at (\x*\spc, 0.00);
  \coordinate (p\x) at (\x*\spc, 0.50);
  \coordinate (f\x) at (\x*\spc,-0.50);
  \coordinate (h\x) at (\x*\spc, 0.875);
}
\node[hnode] (n0) at (c0) {};
\node[dnode] (n1) at (c1) {};
\node[pnode] (n2) at (c2) {};
\node[pnode] (n3) at (c3) {};
\node[dnode] (n4) at (c4) {};
\node[dnode] (n5) at (c5) {};
\node[tnode] (n6) at (c6) {};

\draw[connection] (n0) -- (n1);
\draw[connection] (n1) -- (n2);
\draw[connection] (n2) -- (n3);
\draw[connection] (n3) -- (n4);
\draw[connection] (n4) -- (n5);
\draw[connection] (n5) -- (n6);

\node at (p0) {$(o_2^+)$};
\node at (p1) {$o_2^-$};
\node at (p2) {$o_3^+$};
\node at (p3) {$o_4^+$};
\node at (p4) {$o_4^-$};
\node at (p5) {$o_3^-$};

\node at (f1) {$f_3$};
\node at (f2) {$f_6$};
\node at (f3) {$f_4$};
\node at (f4) {$f_5$};
\node at (f5) {$f_5$};

\draw[couple] (h2) to[out= 45,in= 135] (h5);
\draw[couple] (h3) to[out= 45,in= 135] (h4);

\end{tikzpicture}
\caption{
Example for a route.
The first and the last node is indicated with black and white rectangles, respectively.
Each internal node represents a pickup or delivery,  depicted with black or white circles, respectively.
Above each pickup/delivery node the order, and below it the factory is indicated.}
\label{fig:ex:route}
\end{figure}

Recall that a vehicle's route plan outlines a sequence of factories to visit, along with sets of orders to pickup and deliver at each location, as well as timing information. In this section we use a different, but equivalent representation.
A \emph{route} is a sequence of nodes, where each inner \emph{node} refers to a pickup or a delivery of an order.
That is, a \emph{pickup node} represents the pickup of an order, and a \emph{delivery node} corresponds to the delivery of some order.
The first node of a route is associated with the orders carried by the vehicle, i.e., those orders already picked up but not yet delivered, and the last node marks the end of the route.
We refer to \Cref{fig:ex:route} for an illustration.
It depicts the route of the vehicle corresponding to the center pane of \Cref{fig:ex:states}.
As we can see, the vehicle is on the way to factory~$f_3$ to deliver the carried order~$o_2$.
After that, the vehicle picks up order~$o_3$ at factory~$f_6$ and order~$o_4$ at factory~$f_4$, and then delivers them in LIFO order to factory~$f_5$.

A key gadget of our method is the insertion of an order into a route.
Given a feasible route and some order $o_i$ not in the route, the \emph{insertion} of $o_i$ into the route means that first the pickup node~$o_i^+$ is inserted between two nodes of the route and  then the delivery node~$o_i^-$ is inserted between two nodes of the updated route. 
An insertion is \emph{feasible} if $o_i^+$ precedes $o_i^-$ in the resulting route and the LIFO as well as the capacity constraints are satisfied, and the destination factory does not change (cf. Appendix~\ref{sec:apx:feasibility}).

The \emph{cost} of a solution is calculated by (\ref{eq:obj:manipulated}) throughout this section.
It can computed by the procedure outlined in Appendix~\ref{sec:apx:eval}.

After these preliminaries, our method consists of two main steps:
\begin{enumerate}[1)]
\item Construction of an initial set of routes for the vehicles (Section~\ref{sec:initsol}).
\item Improvement by variable neighborhood search (\Cref{sec:vns}).
\end{enumerate}

\subsubsection{Construction of the initial set of routes}\label{sec:initsol}
Since the solution obtained in the previous epoch could be a good starting point for the current epoch, first we reconstruct and update it.
This involves removing delivery nodes associated with orders already fulfilled and pickups of orders that have already been collected.

Then, new orders are inserted into the updated solution one-by-one.
First, the orders are divided into the sets of urgent and non-urgent orders, respectively. 
The classification is based on the \emph{estimated delay} of an order $o_i \in \tilde{\omega}_k$:
\[
  ed_{k,i} = \left(t^d_i - \tau_k\right) - \left( h^{\operatorname{docking}} + h^p_i + \operatorname{travel}(f^p_i,f^d_i) \right).
\]
The first term represents the remaining time to deliver order~$o_i$ without delay, while the second  term expresses  the minimum time needed to deliver order~$o_i$.
Let $\mathcal{U}_{k} = \{ o_i \in \tilde{\omega}_k : ed_{k,i} \leq U \}$ be the set of urgent orders, and $\overline{\mathcal{U}}_{k} = \tilde{\omega} \setminus \mathcal{U}$  the set of non-urgent ones, where $U$ is a parameter of the algorithm.
 
First, the order in $\mathcal{U}_{k}$ are inserted into the routes  of the vehicles, then those in $\overline{\mathcal{U}}_{k}$.
When processing the orders in either category, orders are grouped by pickup factories, and if two or more orders have the same pickup factory, then priority is given to those order with a larger estimated delay.
For each order $o_i \in \mathcal{U}_{k} \cup \overline{\mathcal{U}}_{k}$ the best feasible insertion is sought. 
That is, $o_i$ is inserted in all feasible ways into the route of each vehicle in turn, and the insertion incurring the least cost increase is chosen. 

\subsubsection{Variable neighborhood search}\label{sec:vns}
In the following, we propose a method based on~VNS to improve the initial set of routes. 
We commence by defining the neighborhood operators to be used within our VNS procedure.

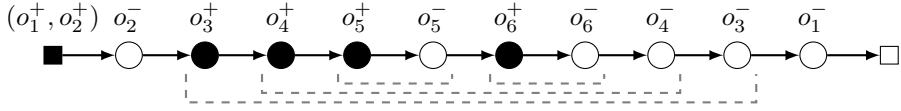
\begin{figure}
\begin{tikzpicture}
\tikzstyle{hnode}=[draw,fill=black,rectangle,minimum height=0.5em, minimum width=0.5em];
\tikzstyle{tnode}=[draw,fill=white,rectangle,minimum height=0.5em, minimum width=0.5em];
\tikzstyle{pnode}=[draw,fill=black,circle,minimum height=1em, minimum width=1em];
\tikzstyle{dnode}=[draw,fill=white,circle,minimum height=1em, minimum width=1em];
\tikzstyle{connection} = [->,thick,latex-,-latex];
\tikzstyle{block} = [thick,dashed,black!50!white];

\def\spc{1.0}

\foreach \x in {0,...,11}
{
  \coordinate (c\x) at (\x*\spc, 0.0);
  \coordinate (p\x) at (\x*\spc, 0.5);
  \coordinate (f\x) at (\x*\spc, -0.5);
}
\node[hnode] (n0) at (c0) {};
\node[dnode] (n1) at (c1) {};
\node[pnode] (n2) at (c2) {};
\node[pnode] (n3) at (c3) {};
\node[pnode] (n4) at (c4) {};
\node[dnode] (n5) at (c5) {};
\node[pnode] (n6) at (c6) {};
\node[dnode] (n7) at (c7) {};
\node[dnode] (n8) at (c8) {};
\node[dnode] (n9) at (c9) {};
\node[dnode] (n10) at (c10) {};
\node[tnode] (n11) at (c11) {};

\draw[connection] (n0) -- (n1);
\draw[connection] (n1) -- (n2);
\draw[connection] (n2) -- (n3);
\draw[connection] (n3) -- (n4);
\draw[connection] (n4) -- (n5);
\draw[connection] (n5) -- (n6);
\draw[connection] (n6) -- (n7);
\draw[connection] (n7) -- (n8);
\draw[connection] (n8) -- (n9);
\draw[connection] (n9) -- (n10);
\draw[connection] (n10) -- (n11);

\node at (p0) {$(o_1^+,o_2^+)$};
\node at (p1) {$o_2^-$};
\node at (p2) {$o_3^+$};
\node at (p3) {$o_4^+$};
\node at (p4) {$o_5^+$};
\node at (p5) {$o_5^-$};
\node at (p6) {$o_6^+$};
\node at (p7) {$o_6^-$};
\node at (p8) {$o_4^-$};
\node at (p9) {$o_3^-$};
\node at (p10) {$o_1^-$};

\draw[block] ($(n4)+(-0.25,-0.25)$) -- ($(n4)+(-0.25,-0.375)$) -- ($(n5)+(+0.25,-0.375)$) -- ($(n5)+(+0.25,-0.25)$);
\draw[block] ($(n6)+(-0.25,-0.25)$) -- ($(n6)+(-0.25,-0.375)$) -- ($(n7)+(+0.25,-0.375)$) -- ($(n7)+(+0.25,-0.25)$);
\draw[block] ($(n3)+(-0.25,-0.25)$) -- ($(n3)+(-0.25,-0.500)$) -- ($(n8)+(+0.25,-0.500)$) -- ($(n8)+(+0.25,-0.25)$);
\draw[block] ($(n2)+(-0.25,-0.25)$) -- ($(n2)+(-0.25,-0.625)$) -- ($(n9)+(+0.25,-0.625)$) -- ($(n9)+(+0.25,-0.25)$);

\end{tikzpicture}
\caption{
A route with blocks  indicated by dashed lines.
}
\label{fig:ex:block}
\end{figure}

The neighborhood operators rely on the notion of blocks and bridges that we define next.
A \emph{block} is a subsequence of consecutive nodes in a vehicle route such that the first and the last node refer to the pickup and the delivery of the same order, respectively.
Refer to  \Cref{fig:ex:block} for an illustration.
Since orders~$o_1$ and~$o_2$ are already loaded on the vehicle, no pickup nodes, and thus no blocks correspond to these orders.

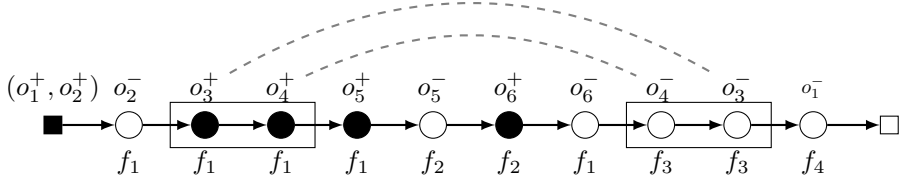
\begin{figure}
\begin{tikzpicture}
\tikzstyle{hnode}=[draw,fill=black,rectangle,minimum height=0.5em, minimum width=0.5em];
\tikzstyle{tnode}=[draw,fill=white,rectangle,minimum height=0.5em, minimum width=0.5em];
\tikzstyle{pnode}=[draw,fill=black,circle,minimum height=1em, minimum width=1em];
\tikzstyle{dnode}=[draw,fill=white,circle,minimum height=1em, minimum width=1em];
\tikzstyle{connection} = [->,thick,latex-,-latex];
\tikzstyle{partner} = [black!50!white];
\tikzstyle{couple} = [thick,dashed,black!50!white];

\def\spc{1.0}

\foreach \x in {0,...,11}
{
  \coordinate (c\x) at (\x*\spc, 0.0);
  \coordinate (p\x) at (\x*\spc, 0.5);
  \coordinate (f\x) at (\x*\spc,-0.5);
}
\node[hnode] (n0) at (c0) {};
\node[dnode] (n1) at (c1) {};
\node[pnode] (n2) at (c2) {};
\node[pnode] (n3) at (c3) {};
\node[pnode] (n4) at (c4) {};
\node[dnode] (n5) at (c5) {};
\node[pnode] (n6) at (c6) {};
\node[dnode] (n7) at (c7) {};
\node[dnode] (n8) at (c8) {};
\node[dnode] (n9) at (c9) {};
\node[dnode] (n10) at (c10) {};
\node[tnode] (n11) at (c11) {};

\node[fit=($(n1)!.66!(n2)$)(n2)(n3)($(n3)!.33!(n4)$),draw] (bleft) {};
\node[fit=($(n7)!.66!(n8)$)(n8)(n9)($(n9)!.33!(n10)$),draw] (bright) {};

\draw[connection] (n0) -- (n1);
\draw[connection] (n1) -- (n2);
\draw[connection] (n2) -- (n3);
\draw[connection] (n3) -- (n4);
\draw[connection] (n4) -- (n5);
\draw[connection] (n5) -- (n6);
\draw[connection] (n6) -- (n7);
\draw[connection] (n7) -- (n8);
\draw[connection] (n8) -- (n9);
\draw[connection] (n9) -- (n10);
\draw[connection] (n10) -- (n11);

\node at (p0) {$(o_1^+,o_2^+)$};
\node at (p1) {$o_2^-$};
\node (c1) at (p2) {$o_3^+$};
\node (c2) at (p3) {$o_4^+$};
\node at (p4) {$o_5^+$};
\node at (p5) {$o_5^-$};
\node at (p6) {$o_6^+$};
\node at (p7) {$o_6^-$};
\node (c3) at (p8) {$o_4^-$};
\node (c4) at (p9) {$o_3^-$};
\node[scale=0.75] at (p10) {$o_1^-$};

\draw[couple] (c1) to[out= 30,in= 150] (c4);
\draw[couple] (c2) to[out= 25,in= 155] (c3);

\node at (f1) {$f_1$};
\node at (f2) {$f_1$};
\node at (f3) {$f_1$};
\node at (f4) {$f_1$};
\node at (f5) {$f_2$};
\node at (f6) {$f_2$};
\node at (f7) {$f_1$};
\node at (f8) {$f_3$};
\node at (f9) {$f_3$};
\node at (f10) {$f_4$};
\end{tikzpicture}
\caption{
Example for a (maximal) bridge.
Left and right sequences of the bridge are framed by rectangles.
}
\label{fig:ex:bridge}
\end{figure}

A \emph{bridge} of a route is  a subsequence of consecutive pickup nodes $(\ell_1,\ldots,\ell_k)$, all belonging the the same factory, and a subsequence  of consecutive delivery nodes $(r_k,\ldots,r_1)$, all delivered to the same factory, where $\ell_j$ and $r_j$ are the pickup and the delivery node of the same order, respectively, for $j=1,\ldots,k$.
An example is depicted in \Cref{fig:ex:bridge}.
A bridge $\left((\ell_1,\ldots,\ell_k)\right.$,$\left.(r_k,\ldots,r_1)\right)$ is \emph{maximal}, if there is no bridge containing it properly.
That is, neither
$\left((\operatorname{pred}(\ell_1),\ell_1,\ldots,\ell_k)\right.$, $\left.(r_k,\ldots,r_1,\operatorname{succ}(r_1)\right)$
nor
$\left((\ell_1,\ldots,\ell_k,\operatorname{succ}(\ell_k)\right.$, $\left.(\operatorname{pred}(r_k),r_k,\ldots,r_1)\right)$
constitutes a bridge, where $\operatorname{pred}(n)$ and $\operatorname{succ}(n)$ denote the immediate predecessor and the immediate successor of node~$n$ in the route, respectively.

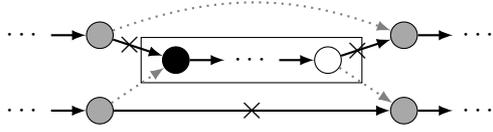
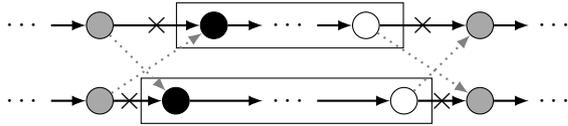
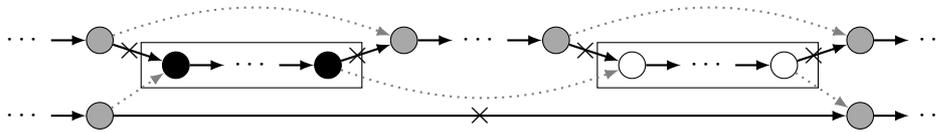
\begin{figure}
\begin{subfigure}[b]{\textwidth}
\begin{tikzpicture}
\tikzstyle{pnode}=[draw,fill=black,circle,minimum height=1em, minimum width=1em];
\tikzstyle{dnode}=[draw,fill=white,circle,minimum height=1em, minimum width=1em];
\tikzstyle{xnode}=[draw,fill=white!66!black,circle,minimum height=1em, minimum width=1em];
\tikzstyle{connection} = [->,thick,latex-,-latex];
\tikzstyle{newconnection} = [connection,dotted,white!50!black];

\node        (a1) at (-2,0) {$\cdots$};
\node[xnode] (a2) at (-1,0) {};
\node[pnode] (a3) at ( 0,-0.33) {};
\node        (a4) at ( 1,-0.33) {$\cdots$};
\node[dnode] (a5) at ( 2,-0.33) {};
\node[xnode] (a6) at ( 3,0) {};
\node        (a7) at ( 4,0) {$\cdots$};

\draw[connection] (a1) -- (a2);
\draw[connection] (a2) -- (a3) node[pos=0.33,scale=1.25]{$\times$};
\draw[connection] (a3) -- (a4);
\draw[connection] (a4) -- (a5);
\draw[connection] (a5) -- (a6) node[pos=0.33,scale=1.25]{$\times$};
\draw[connection] (a6) -- (a7);

\node[fit=($(a2)!.66!(a3)$)(a3)(a5)($(a5)!.33!(a6)$),draw] {};

\node        (b1) at (-2,-1) {$\cdots$};
\node[xnode] (b2) at (-1,-1) {};
\node[xnode] (b3) at ( 3,-1) {};
\node        (b4) at ( 4,-1) {$\cdots$};

\draw[connection] (b1) -- (b2);
\draw[connection] (b2) -- (b3) node[pos=0.5,scale=1.25]{$\times$};
\draw[connection] (b3) -- (b4);

\draw[newconnection] (a2) to [out=20,in=160] (a6);
\draw[newconnection] (b2) -- (a3);
\draw[newconnection] (a5) -- (b3);
\end{tikzpicture}
\caption{Block relocation between routes.}
\label{fig:ex:ops:blrel}
\end{subfigure}
\par\bigskip
\begin{subfigure}[b]{\textwidth}
\begin{tikzpicture}
\tikzstyle{pnode}=[draw,fill=black,circle,minimum height=1em, minimum width=1em];
\tikzstyle{dnode}=[draw,fill=white,circle,minimum height=1em, minimum width=1em];
\tikzstyle{xnode}=[draw,fill=white!66!black,circle,minimum height=1em, minimum width=1em];
\tikzstyle{connection} = [->,thick,latex-,-latex];
\tikzstyle{newconnection} = [connection,dotted,white!50!black];

\node        (a1) at (-2,0) {$\cdots$};
\node[xnode] (a2) at (-1,0) {};
\node[pnode] (a3) at ( 0.5,0) {};
\node        (a4) at ( 1.5,0) {$\cdots$};
\node[dnode] (a5) at ( 2.5,0) {};
\node[xnode] (a6) at ( 4,0) {};
\node        (a7) at ( 5,0) {$\cdots$};

\draw[connection] (a1) -- (a2);
\draw[connection] (a2) -- (a3) node[pos=0.5,scale=1.25]{$\times$};
\draw[connection] (a3) -- (a4);
\draw[connection] (a4) -- (a5);
\draw[connection] (a5) -- (a6) node[pos=0.5,scale=1.25]{$\times$};
\draw[connection] (a6) -- (a7);

\node        (b1) at (-2,-1) {$\cdots$};
\node[xnode] (b2) at (-1,-1) {};
\node[pnode] (b3) at ( 0,-1) {};
\node        (b4) at (1.5,-1) {$\cdots$};
\node[dnode] (b5) at ( 3,-1) {};
\node[xnode] (b6) at ( 4,-1) {};
\node        (b7) at ( 5,-1) {$\cdots$};

\node[fit=($(a2)!.75!(a3)$)(a3)(a5)($(a5)!.25!(a6)$),draw] {};
\node[fit=($(b2)!.66!(b3)$)(b3)(b5)($(b5)!.25!(b6)$),draw] {};

\draw[connection] (b1) -- (b2);
\draw[connection] (b2) -- (b3) node[pos=0.33,scale=1.25]{$\times$};
\draw[connection] (b3) -- (b4);
\draw[connection] (b4) -- (b5);
\draw[connection] (b5) -- (b6) node[pos=0.5,scale=1.25]{$\times$};
\draw[connection] (b6) -- (b7);

\draw[newconnection] (b2) -- (a3);
\draw[newconnection] (b5) -- (a6);
\draw[newconnection] (a2) -- (b3);
\draw[newconnection] (a5) -- (b6);
\end{tikzpicture}
\caption{Block exchange between routes.}
\label{fig:ex:ops:blexc}
\end{subfigure}
\par\bigskip
\begin{subfigure}[b]{\textwidth}
\begin{tikzpicture}
\tikzstyle{pnode}=[draw,fill=black,circle,minimum height=1em, minimum width=1em];
\tikzstyle{dnode}=[draw,fill=white,circle,minimum height=1em, minimum width=1em];
\tikzstyle{xnode}=[draw,fill=white!66!black,circle,minimum height=1em, minimum width=1em];
\tikzstyle{connection} = [->,thick,latex-,-latex];
\tikzstyle{newconnection} = [connection,dotted,white!50!black];

\node        (a1) at (-2,0) {$\cdots$};
\node[xnode] (a2) at (-1,0) {};
\node[pnode] (a3) at ( 0,-0.33) {};
\node        (a4) at ( 1,-0.33) {$\cdots$};
\node[pnode] (a5) at ( 2,-0.33) {};
\node[xnode] (a6) at ( 3,0) {};
\node        (a7) at ( 4,0) {$\cdots$};
\node[xnode] (a8) at (5,0) {};
\node[dnode] (a9) at (6,-0.33) {};
\node        (a10) at (7,-0.33) {$\cdots$};
\node[dnode] (a11) at (8,-0.33) {};
\node[xnode] (a12) at (9,0) {};
\node        (a13) at (10,0) {$\cdots$};

\draw[connection] (a1) -- (a2);
\draw[connection] (a2) -- (a3) node[pos=0.33,scale=1.25]{$\times$};
\draw[connection] (a3) -- (a4);
\draw[connection] (a4) -- (a5);
\draw[connection] (a5) -- (a6) node[pos=0.33,scale=1.25]{$\times$};
\draw[connection] (a6) -- (a7);
\draw[connection] (a7) -- (a8);
\draw[connection] (a8) -- (a9) node[pos=0.33,scale=1.25]{$\times$};
\draw[connection] (a9) -- (a10);
\draw[connection] (a10) -- (a11);
\draw[connection] (a11) -- (a12) node[pos=0.33,scale=1.25]{$\times$};
\draw[connection] (a12) -- (a13);

\node[fit=($(a2)!.66!(a3)$)(a3)(a5)($(a5)!.33!(a6)$),draw] {};
\node[fit=($(a8)!.66!(a9)$)(a9)(a11)($(a11)!.33!(a12)$),draw] {};

\node        (b1) at (-2,-1) {$\cdots$};
\node[xnode] (b2) at (-1,-1) {};
\node[xnode] (b6) at (9,-1) {};
\node        (b7) at (10,-1) {$\cdots$};

\draw[connection] (b1) -- (b2);
\draw[connection] (b6) -- (b7);

\draw[connection] (b2) -- (b6) node[pos=0.5,scale=1.25]{$\times$};

\draw[newconnection] (a2) to [out=20,in=160] (a6);
\draw[newconnection] (b2) -- (a3);

\draw[newconnection] (a8) to [out=20,in=160] (a12);
\draw[newconnection] (a11) -- (b6);

\draw[newconnection] (a5) to [out=-20,in=-160] (a9);
\end{tikzpicture}
\caption{
Bridge relocation between routes.
The two parts of the bridge will be adjacent after relocation.
}
\label{fig:ex:ops:brrel}
\end{subfigure}
\caption{
Example for schedule operations.
Crossed/dotted arrows refer to old/new links.   
}
\label{fig:ex:ops}
\end{figure}

A \emph{neighborhood operator} perturbs a solution slightly and provides a modified solution.
The following neighborhood operators will be applied to the solutions:
\begin{enumerate}[i)]
\item The \emph{relocate-block} operator selects a vehicle route and a block in the route, and relocates it to another position in the same or in another route (cf.~\Cref{fig:ex:ops:blrel}).

\item The \emph{block-exchange} operator selects two disjoint blocks from the same or from distinct vehicle routes, and exchanges them as depicted in \Cref{fig:ex:ops:blexc}.

\item The \emph{relocate-bridge} operator selects a maximal bridge from a vehicle route and moves it into another position of the same route or to a different route as shown in  \Cref{fig:ex:ops:brrel}.
\end{enumerate}

When applying an operator, the new position for a block or bridge is always chosen in such a manner that the resulting solution satisfies the LIFO constraint.
However, if any other constraint is violated, the solution is dropped and the operator fails.

\begin{algorithm}
\caption{VNS procedure}\label{algo:vns}
\begin{algorithmic}[1]
\While{termination condition is not met}\label{vns:test}
\State Find the best neighbor $S'$ of $S$ using the 
{\em relocate-bridge} operator.
\If{ $S'$ has a lower cost than $S$}
\State $S := S'$ and proceed with Step~\ref{vns:test}.
\EndIf
\State Find the best neighbor $S'$ of $S$ using the 
{\em block-exchange} operator.
\If{ $S'$ has a lower cost than $S$}
\State $S := S'$ and proceed with Step~\ref{vns:test}.
\EndIf
\State Find the best neighbor $S'$ of $S$ using the 
{\em relocate-block} operator.
\If{ $S'$ has a lower cost than $S$}
\State $S := S'$ and proceed with Step~\ref{vns:test}.
\Else
\State Proceed with Step~\ref{vns:output}.
\EndIf
\EndWhile
\State Output $S$.\label{vns:output}
\end{algorithmic}
\end{algorithm}

Our VNS procedure, depicted in \Cref{algo:vns}, receives an initial feasible solution and tries to improve it repeatedly by applying the neighborhood operators.
Throughout the algorithm, $S$ denotes the current solution.
The termination condition can be a time limit, or the number of iterations of the procedure.
Solutions are compared based on their cost computed by the reward function (\ref{eq:obj:manipulated}).
The neighborhood operators are applied exhaustively, i.e., they are applied in all possible ways to find the best neighboring solution. For instance, in case of the relocate-bridge operator, when seeking the best neighbor of~$S$, in each route all maximal bridges are identified, and they are reinserted in the same or in a different route in all possible ways, and then the insertion producing the least cost feasible solution is chosen.
If an operator cannot be applied, or fails to produce any feasible neighbor, then we assume that~$S'$ has a larger cost than~$S$ and proceed accordingly.
Since only improvement of the current solution is allowed, cycling cannot occur, and the procedure terminates in a local minimum.

\section{Computational experiments}\label{sec:comp}
In this section, we present the results our computational experiments.
In \Cref{sec:pre:dpdp}, we describe our case study.
In \Cref{sec:comp:instances}, we introduce the benchmark dataset.
In \Cref{sec:comp:eval}, we tune the objective multipliers to get an efficient solution approach.
Finally, we compare our best approach to other methods in \Cref{sec:comp:comp}.

\paragraph{Setup}
All our experiments were performed on a workstation with an Intel Core i9-7960X 2.80~GHz CPU with 16~cores, under Debian~9 operating system using a single thread.
Due to the 4~hours time-window for the orders (see \Cref{sec:comp:instances}), we used parameter $U=3600$ when urgent and non-urgent orders are determined for initial solution.
Since each epoch is 10~minutes long, we applied a time limit of 9~minutes in our variable neighborhood search.

\subsection{Case study: The Dynamic Pickup and Delivery Problem challenge}\label{sec:pre:dpdp}
We consider \emph{The Dynamic Pickup and Delivery Problem challenge} as a case study \citep{hao2022introduction}.
This challenge was organized as a competition of the International Conference on Automated Planning and Scheduling (ICAPS) in~2021 with 152 participating teams.
The best three teams (i.e., the teams whose algorithms achieved the smallest scores on a hidden dataset) presented their solution approaches in the conference and the challenge generated a series of papers in the scientific literature in the recent years.

The gold medal team, \citet{goldsolution} proposed an RO method for the problem, where a variable neighborhood search was used to optimize the states with the currently known orders.
In each epoch, the authors reconstructed the solution from the previous period, and dispatched new orders with a cheapest insertion heuristics.
Four local search operators (couple-exchange, block-exchange, relocate-block, and multi-relocate) were used in a variable neighborhood search together with a disturbance operator to perturb solutions.
Since then, the team published their results in \citep{cai2022variable}.
The silver medal team, \citet{silversolution} developed an RO method for the problem, where a rule-based procedure was applied to dispatch orders.
The bronze-medal team, \citet{bronzesolution} proposed a CFA method for the problem, where a local search-based approach was used to solve the problem with a perturbed objective function.
In each epoch, the authors inserted non-dispatched orders into the solution with a cheapest insertion heuristic.
Three local search operators (block-exchange, relocate-block, and a custom made) were used to improve solutions.
An extra cost term was introduced to penalize if a factory is visited by too many vehicles, thereby reducing the chance of waiting for free docking ports.

\citet{cai2022efficient} propose an RO method for the~DPDP, where a reference point-based multi-objective evolutionary algorithm is used to solve the states.
Four local search operators are used in a variable neighborhood search to improve solutions as in \citep{cai2022variable}.
The authors compare their solution approach to the algorithms of the podium teams.
For the comparison, the authors use the benchmark dataset with the exception of the largest instances.
Recently, \citet{cai2023survey} provide a review of the dynamic pickup and delivery problem literature covering the last two decades.
As a case study, the authors provide the comparison of \citet{cai2022efficient} on the full benchmark dataset.

\citet{du2023hierarchical} propose an RO method for the~DPDP with a hierarchical optimization approach.
The authors are the first to model the problem as an~SDP.
They apply several order dispatching strategies to assign orders to vehicles.
A buffering pool is also used to postpone the assignment of some non-urgent requests.
A local search operator is introduced to improve the solution.
The authors compare their solution approach to a greedy baseline strategy, and to the  silver-winning algorithm of \citet{silversolution} (which was falsely claimed to be the best algorithm of the competition).
For the comparison, the authors use only a selection of the benchmark dataset, but also generated new problem instances.

\subsection{Instances}\label{sec:comp:instances}

\begin{table}
\caption{Basic properties of the benchmark dataset.}
\label{tab:dataset:public}
\begin{threeparttable}
\begin{tabular}{rrrr}
\toprule
Group & Instances & Orders & Vehicles\\
\midrule
1 &  1 --  8 &   50 &   5 \\
2 &  9 -- 16 &  100 &   5 \\
3 & 17 -- 24 &  300 &  20 \\
4 & 25 -- 32 &  500 &  20 \\
5 & 33 -- 40 & 1000 &  50 \\
6 & 41 -- 48 & 2000 &  50 \\
7 & 49 -- 56 & 3000 & 100 \\
8 & 57 -- 64 & 4000 & 100 \\
\bottomrule
\end{tabular}
\end{threeparttable}
\end{table}

We use the publicly available dataset of the ICAPS 2021 DPDP competition \citep{hao2022introduction}.
This dataset consists of 64 instances based on 30~days of historical data of Huawei.
These instances contain 50-4000 orders of a single day to be satisfied with 5-100 vehicles, see \Cref{tab:dataset:public}.
In the following, we briefly describe the public instances, while for detailed description we refer to \citep{hao2022introduction}.

In case of all instances, there are 4~hours to complete an order on time, that is, for each order~$o_i\in\mathcal{O}$ we have $t^d_i - t^p_i = 14\,400$ seconds.
Each order is given as a set of \emph{order items}, where each item is either a \emph{box}, a \emph{small pallet}, or a \emph{standard pallet}, with quantity of 0.25, 0.5, and 1~unit, respectively, while the uniform capacity limit of the vehicles is 15~units.
Loading or unloading a box, a small pallet, and a standard pallet takes 15, 30, and 60~seconds, respectively.
The underlying network is the same for all instances and consists of 153~factories, where each factory has 6~docking ports.
Docking to a port takes  half an hour (i.e., $h^{\operatorname{docking}} = 1800$ seconds).
Finally, $\lambda_1 = 1/n$ and $\lambda_2 = 10\,000/3600 \approx 2.78$ in the objective function~\eqref{eq:obj}, that is, a cost of~10000 monetary units must be paid for every hour of delay.

Note that the total size of the  items of an order may exceed the uniform capacity of the vehicles.
In this and only this case, the order items can be transported separately.
The completion time of such an order is the latest completion time of its order items.
For such big orders, we arranged items in a non-increasing size order, then we applied a first fit procedure to divide items into separate orders.
By this, we got the same problem as proposed earlier, with the tiny difference that the tardiness has to be calculated differently for split orders.

\paragraph{Dynamic evaluation}
The organizers of the ICAPS 2021 DPDP competition also provided a simulator to support the dynamic evaluation of the solution approaches.
This simulator essentially follows the sequential decision procedure described in \Cref{sec:SDP}.
The operating horizon is divided into epochs of length 10 minutes each.
The final output of the simulator is the score of the dispatching algorithm on the given instance.
For details, we refer to \citep{hao2022introduction}.

\subsection{Evaluation of the cost function approximation method}\label{sec:comp:eval}
In the following, we evaluate our method with different $\lambda_3$ and $\lambda_4$ parameters with the aim of finding the best settings.

\subsubsection{Significance of penalizing the waiting times}
In these experiments, we suppressed the fourth term in the objective function~(\ref{eq:obj:manipulated}).
That is, we ran our method with objective function coefficients $\lambda_1 = 1/n$, $\lambda_2 = 2.78$, $\lambda_4 = 0$, and $\lambda_3$ chosen from $\{0.0, 0.25\times \lambda_2,0.5\times \lambda_2, 0.75\times\lambda_2\}$. 
The average scores for each group, and for all instances are indicated in \Cref{tab:eval:waiting}, while detailed results can be found in the supplementary data.

\begin{table}
\caption{Evaluation of penalizing waiting times.}
\label{tab:eval:waiting}
\begin{threeparttable}
\begin{tabular}{r rrrr}
\toprule
& \multicolumn{4}{c}{Multiplier for waiting times ($\lambda_3$)}\\
\cmidrule(lr){2-5}
Group & $0\times\lambda_2$ & $0.25\times\lambda_2$ & $0.5\times\lambda_2$ & $0.75\times\lambda_2$\\
\midrule
1 &       1\,307.5 &      1\,307.5 &      1\,307.5 &      1\,307.5\\
2 &      31\,707.6 &     31\,707.6 &      31 707.6 &      31 707.6\\
\addlinespace[1ex]
3 &          690.0 &         690.0 &         690.0 &         690.0\\
4 &       6\,750.9 &      6\,750.5 &      6\,750.5 &      6\,750.5\\
\addlinespace[1ex]
5 &      11\,597.1 &     11\,225.6 &     11\,225.6 &     11\,225.6\\
6 &      49\,011.9 &     53\,183.2 &     54\,091.5 &     49\,841.4\\
\addlinespace[1ex]
7 &  1\,076\,526.1 &    590\,379.6 &    738\,200.9 &    769\,756.1\\
8 & 11\,605\,754.6 & 4\,696\,315.2 & 4\,618\,058.6 & 4\,412\,835.5\\
\cmidrule(lr){1-5}
average: & 1\,597\,918.2 & 673\,944.9 & 682\,754.0 & 660\,514.3\\
\bottomrule
\end{tabular}
\end{threeparttable}
\end{table}

Clearly, penalizing waiting times has no impact on the first two groups, since the number of docking ports is greater than the number of vehicles, and thus no waiting occurs.
On instances with 20 and 50 vehicles (Groups 3-6), the differences are negligible.
However, in case of 100 vehicles (Groups 7-8), penalizing idle times significantly reduced the average score.
On the largest group, the average improvement is 60-62\%.

\begin{figure}
\centering
\begin{subfigure}[b]{0.45\textwidth}
\begin{tikzpicture}
\begin{axis}[
    width=\textwidth,
    height=6cm,
    ymajorgrids,
    xmin=0,
    xmax=200,
    ymin=0,
    ymax=100,
    yticklabel pos=left,
    ytick = {10,33,50,100},
    yticklabel={$\pgfmathprintnumber{\tick}$},
    xlabel=Update points,
    ylabel=Number of waiting vehicles,
    legend pos= north west
]

\addplot[red] coordinates {
 (  0, 0) (  1, 0) (  2, 0) (  3, 0) (  4, 0) (  5, 0) (  6, 0) (  7, 0) (  8, 0) (  9, 0)
 ( 10, 0) ( 11, 0) ( 12, 0) ( 13, 0) ( 14, 0) ( 15, 0) ( 16, 0) ( 17, 0) ( 18, 0) ( 19, 1)
 ( 20, 0) ( 21, 0) ( 22, 0) ( 23, 0) ( 24, 0) ( 25, 1) ( 26, 1) ( 27, 0) ( 28, 0) ( 29, 2)
 ( 30, 1) ( 31, 1) ( 32, 1) ( 33, 2) ( 34, 2) ( 35, 2) ( 36, 0) ( 37, 0) ( 38, 0) ( 39, 0)
 ( 40, 3) ( 41, 3) ( 42, 3) ( 43, 3) ( 44, 3) ( 45, 5) ( 46, 5) ( 47, 2) ( 48, 4) ( 49, 6)
 ( 50, 6) ( 51, 6) ( 52, 6) ( 53, 8) ( 54,10) ( 55, 9) ( 56, 6) ( 57, 9) ( 58,10) ( 59,12)
 ( 60,10) ( 61, 9) ( 62,13) ( 63,11) ( 64, 9) ( 65, 7) ( 66, 7) ( 67, 2) ( 68, 4) ( 69, 6)
 ( 70, 5) ( 71, 3) ( 72, 3) ( 73, 4) ( 74, 3) ( 75, 1) ( 76, 2) ( 77, 0) ( 78, 0) ( 79, 1)
 ( 80, 0) ( 81, 0) ( 82, 1) ( 83, 5) ( 84, 4) ( 85, 1) ( 86, 1) ( 87, 6) ( 88, 6) ( 89, 6)
 ( 90, 7) ( 91,10) ( 92, 7) ( 93,12) ( 94,12) ( 95,14) ( 96,12) ( 97,11) ( 98,11) ( 99,12)
 (100, 9) (101,10) (102, 9) (103, 9) (104, 9) (105, 9) (106, 8) (107,11) (108,12) (109,15)
 (110,16) (111,18) (112,21) (113,22) (114,20) (115,23) (116,27) (117,27) (118,28) (119,24)
 (120,31) (121,31) (122,36) (123,34) (124,36) (125,36) (126,42) (127,48) (128,47) (129,46)
 (130,46) (131,42) (132,43) (133,45) (134,44) (135,40) (136,41) (137,40) (138,45) (139,45)
 (140,45) (141,43) (142,47) (143,44) (144,43) (145,45) (146,45) (147,40) (148,43) (149,42)
 (150,42) (151,37) (152,38) (153,37) (154,39) (155,42) (156,44) (157,43) (158,44) (159,40)
 (160,46) (161,43) (162,46) (163,44) (164,38) (165,37) (166,34) (167,34) (168,31) (169,29)
 (170,27) (171,27) (172,23) (173,20) (174,17) (175,17) (176,16) (177,14) (178,13) (179,11)
 (180,11) (181, 8) (182, 6) (183, 5) (184, 3) (185, 2) (186, 1) (187, 0) (188, 0) (189, 0)
 (190, 0) (191, 0)
};

\addplot[blue] coordinates {
 (  0, 0) (  1, 0) (  2, 0) (  3, 0) (  4, 0) (  5, 0) (  6, 0) (  7, 0) (  8, 0) (  9, 0)
 ( 10, 0) ( 11, 0) ( 12, 0) ( 13, 0) ( 14, 0) ( 15, 0) ( 16, 0) ( 17, 0) ( 18, 0) ( 19, 0)
 ( 20, 0) ( 21, 0) ( 22, 0) ( 23, 0) ( 24, 0) ( 25, 0) ( 26, 0) ( 27, 0) ( 28, 0) ( 29, 0)
 ( 30, 0) ( 31, 0) ( 32, 0) ( 33, 0) ( 34, 0) ( 35, 0) ( 36, 0) ( 37, 0) ( 38, 0) ( 39, 0)
 ( 40, 0) ( 41, 1) ( 42, 0) ( 43, 0) ( 44, 1) ( 45, 0) ( 46, 0) ( 47, 0) ( 48, 0) ( 49, 1)
 ( 50, 1) ( 51, 1) ( 52, 0) ( 53, 1) ( 54, 0) ( 55, 1) ( 56, 0) ( 57, 0) ( 58, 1) ( 59, 1)
 ( 60, 2) ( 61, 1) ( 62, 1) ( 63, 0) ( 64, 1) ( 65, 0) ( 66, 0) ( 67, 0) ( 68, 0) ( 69, 1)
 ( 70, 0) ( 71, 1) ( 72, 0) ( 73, 2) ( 74, 1) ( 75, 0) ( 76, 0) ( 77, 0) ( 78, 1) ( 79, 0)
 ( 80, 1) ( 81, 0) ( 82, 0) ( 83, 0) ( 84, 1) ( 85, 2) ( 86, 1) ( 87, 0) ( 88, 0) ( 89, 0)
 ( 90, 3) ( 91, 1) ( 92, 3) ( 93, 3) ( 94, 1) ( 95, 1) ( 96, 3) ( 97, 1) ( 98, 3) ( 99, 3)
 (100, 4) (101, 5) (102, 2) (103, 1) (104, 3) (105, 1) (106, 0) (107, 0) (108, 1) (109, 2)
 (110, 0) (111, 1) (112, 6) (113, 4) (114, 5) (115, 4) (116, 0) (117, 3) (118, 3) (119, 2)
 (120, 4) (121, 3) (122, 3) (123, 7) (124, 9) (125, 9) (126, 6) (127,10) (128, 8) (129, 7)
 (130, 8) (131, 8) (132, 9) (133,10) (134, 9) (135,10) (136, 8) (137, 8) (138, 7) (139, 5)
 (140, 4) (141, 4) (142, 5) (143, 5) (144, 6) (145, 7) (146, 4) (147, 5) (148, 4) (149, 4)
 (150, 5) (151, 3) (152, 2) (153, 1) (154, 1) (155, 1) (156, 1) (157, 1) (158, 3) (159, 2)
 (160, 1) (161, 0) (162, 0) (163, 1) (164, 0) (165, 0) (166, 2) (167, 0) (168, 0) (169, 0)
 (170, 0) (171, 0) (172, 0) (173, 0) (174, 0) (175, 0) (176, 0) (177, 0) (178, 0) (179, 0)
 (180, 0) (181, 0) (182, 0) (183, 0) (184, 0)
};

\legend{RO,CFA}

\end{axis}
\end{tikzpicture}
\caption{Number of waiting vehicles.}
\label{fig:waiting:vehicles}
\end{subfigure}
\hfill
\begin{subfigure}[b]{0.45\textwidth}
\begin{tikzpicture}
\begin{axis}[
    width=\textwidth,
    height=6cm,
    ymajorgrids,
    xmin=0,
    xmax=200,
    ymin=0,
    ymax=120,
    yticklabel pos=left,
    ytick= {30,45,60,90,120},
    yticklabel={$\pgfmathprintnumber{\tick}$\,min},
    xlabel=Update points,
    ylabel=Average committed time,
    legend pos= north west
]

\addplot[red] coordinates {
 (  0,   0.0) (  1,   0.0) (  2,   2.4) (  3,   3.0) (  4,   4.8) (  5,   7.0) (  6,   7.2) (  7,   8.3) (  8,   7.4) (  9,   7.4)
 ( 10,   7.5) ( 11,   8.3) ( 12,   7.5) ( 13,   8.2) ( 14,  11.3) ( 15,  10.1) ( 16,  10.2) ( 17,  10.1) ( 18,  13.1) ( 19,  13.7)
 ( 20,  14.8) ( 21,  14.9) ( 22,  18.5) ( 23,  16.1) ( 24,  15.6) ( 25,  18.5) ( 26,  20.4) ( 27,  25.1) ( 28,  24.0) ( 29,  22.4)
 ( 30,  25.6) ( 31,  29.7) ( 32,  26.5) ( 33,  23.0) ( 34,  26.3) ( 35,  29.7) ( 36,  29.6) ( 37,  29.6) ( 38,  28.9) ( 39,  28.4)
 ( 40,  27.6) ( 41,  29.2) ( 42,  27.9) ( 43,  26.8) ( 44,  30.4) ( 45,  30.7) ( 46,  31.9) ( 47,  29.7) ( 48,  31.5) ( 49,  30.2)
 ( 50,  31.1) ( 51,  32.4) ( 52,  33.2) ( 53,  36.6) ( 54,  32.7) ( 55,  33.8) ( 56,  32.1) ( 57,  29.4) ( 58,  27.1) ( 59,  26.5)
 ( 60,  27.9) ( 61,  33.7) ( 62,  31.2) ( 63,  34.1) ( 64,  31.6) ( 65,  32.4) ( 66,  30.9) ( 67,  30.5) ( 68,  28.0) ( 69,  27.3)
 ( 70,  27.5) ( 71,  31.9) ( 72,  31.1) ( 73,  29.2) ( 74,  25.3) ( 75,  29.0) ( 76,  28.6) ( 77,  26.1) ( 78,  26.5) ( 79,  27.5)
 ( 80,  34.9) ( 81,  34.8) ( 82,  33.2) ( 83,  34.9) ( 84,  32.4) ( 85,  30.8) ( 86,  28.1) ( 87,  31.1) ( 88,  31.0) ( 89,  36.8)
 ( 90,  40.0) ( 91,  42.1) ( 92,  40.1) ( 93,  35.7) ( 94,  35.7) ( 95,  38.3) ( 96,  37.1) ( 97,  35.9) ( 98,  33.8) ( 99,  38.6)
 (100,  39.3) (101,  45.8) (102,  45.3) (103,  46.6) (104,  46.9) (105,  46.1) (106,  43.8) (107,  44.9) (108,  54.9) (109,  55.0)
 (110,  51.4) (111,  57.9) (112,  56.5) (113,  64.5) (114,  61.4) (115,  67.7) (116,  70.5) (117,  72.3) (118,  70.8) (119,  75.8)
 (120,  78.1) (121,  91.0) (122,  83.5) (123,  85.4) (124,  86.4) (125,  85.4) (126,  92.9) (127,  95.1) (128,  91.4) (129,  87.9)
 (130,  88.9) (131,  92.9) (132,  95.1) (133,  88.9) (134,  85.8) (135,  88.8) (136,  82.7) (137,  93.2) (138,  98.5) (139, 109.5)
 (140, 107.2) (141, 107.8) (142, 110.6) (143, 108.6) (144, 110.4) (145, 106.1) (146, 104.4) (147, 102.7) (148,  98.2) (149,  96.6)
 (150,  89.7) (151,  91.7) (152,  95.0) (153,  97.6) (154,  95.4) (155,  92.2) (156,  88.5) (157,  86.0) (158,  83.2) (159,  87.1)
 (160,  81.2) (161,  78.8) (162,  75.8) (163,  71.1) (164,  67.8) (165,  66.0) (166,  65.4) (167,  63.2) (168,  58.9) (169,  53.1)
 (170,  46.2) (171,  40.8) (172,  35.6) (173,  32.8) (174,  30.3) (175,  27.0) (176,  23.6) (177,  20.3) (178,  17.4) (179,  14.3)
 (180,  11.5) (181,  10.7) (182,   9.4) (183,   7.1) (184,   5.2) (185,   3.5) (186,   2.2) (187,   1.3) (188,   0.6) (189,   0.2)
 (190,   0.0) (191,   0.0)
};

\addplot[blue] coordinates {
 (  0,   0.0) (  1,   0.0) (  2,   2.4) (  3,   3.0) (  4,   4.8) (  5,   7.0) (  6,   7.2) (  7,   8.3) (  8,   7.4) (  9,   7.4)
 ( 10,   7.5) ( 11,   8.3) ( 12,   7.5) ( 13,   8.2) ( 14,  11.3) ( 15,  10.1) ( 16,  10.2) ( 17,  10.1) ( 18,  12.2) ( 19,  12.9)
 ( 20,  14.9) ( 21,  14.9) ( 22,  18.2) ( 23,  17.7) ( 24,  16.1) ( 25,  18.2) ( 26,  16.6) ( 27,  22.2) ( 28,  21.9) ( 29,  20.0)
 ( 30,  23.6) ( 31,  24.3) ( 32,  26.3) ( 33,  23.1) ( 34,  26.0) ( 35,  29.7) ( 36,  27.5) ( 37,  28.3) ( 38,  27.5) ( 39,  24.6)
 ( 40,  28.9) ( 41,  27.5) ( 42,  24.6) ( 43,  26.5) ( 44,  24.2) ( 45,  28.6) ( 46,  28.5) ( 47,  27.1) ( 48,  27.2) ( 49,  25.9)
 ( 50,  27.1) ( 51,  31.1) ( 52,  33.0) ( 53,  30.8) ( 54,  27.5) ( 55,  29.3) ( 56,  30.6) ( 57,  29.0) ( 58,  26.4) ( 59,  23.6)
 ( 60,  21.3) ( 61,  29.3) ( 62,  27.8) ( 63,  29.6) ( 64,  29.5) ( 65,  29.5) ( 66,  28.5) ( 67,  27.8) ( 68,  29.2) ( 69,  28.4)
 ( 70,  30.2) ( 71,  33.0) ( 72,  32.0) ( 73,  31.1) ( 74,  31.8) ( 75,  30.7) ( 76,  28.2) ( 77,  25.9) ( 78,  26.2) ( 79,  29.7)
 ( 80,  34.3) ( 81,  33.8) ( 82,  33.1) ( 83,  33.0) ( 84,  31.3) ( 85,  30.2) ( 86,  28.6) ( 87,  31.8) ( 88,  29.8) ( 89,  30.9)
 ( 90,  29.8) ( 91,  30.9) ( 92,  26.6) ( 93,  27.3) ( 94,  27.9) ( 95,  30.6) ( 96,  27.0) ( 97,  31.8) ( 98,  32.6) ( 99,  32.0)
 (100,  35.8) (101,  34.2) (102,  36.3) (103,  35.2) (104,  34.6) (105,  35.4) (106,  37.5) (107,  36.6) (108,  38.8) (109,  35.0)
 (110,  41.4) (111,  38.3) (112,  36.3) (113,  39.1) (114,  36.8) (115,  40.0) (116,  43.0) (117,  43.1) (118,  43.4) (119,  41.6)
 (120,  44.4) (121,  45.4) (122,  44.4) (123,  38.9) (124,  40.5) (125,  39.8) (126,  43.0) (127,  42.3) (128,  45.5) (129,  42.4)
 (130,  42.7) (131,  37.4) (132,  37.7) (133,  45.9) (134,  42.1) (135,  37.5) (136,  34.9) (137,  33.0) (138,  38.0) (139,  37.4)
 (140,  35.8) (141,  34.6) (142,  34.8) (143,  35.5) (144,  34.8) (145,  35.6) (146,  34.4) (147,  31.6) (148,  34.3) (149,  30.9)
 (150,  25.4) (151,  30.7) (152,  30.8) (153,  29.3) (154,  29.8) (155,  30.3) (156,  27.6) (157,  26.0) (158,  24.9) (159,  25.0)
 (160,  24.8) (161,  26.9) (162,  25.2) (163,  21.3) (164,  19.4) (165,  18.8) (166,  16.1) (167,  15.7) (168,  12.7) (169,  10.2)
 (170,   7.5) (171,   5.1) (172,   4.0) (173,   3.6) (174,   3.3) (175,   2.5) (176,   2.0) (177,   1.5) (178,   1.1) (179,   0.8)
 (180,   0.6) (181,   0.4) (182,   0.2) (183,   0.1) (184,   0.0)
};

\legend{RO,CFA}

\end{axis}
\end{tikzpicture}
\caption{Average committed time.}
\label{fig:waiting:time}
\end{subfigure}
\caption{
Results on a large instance.
Red and blue plots refer to the case where decision problems were solved with the original reward function~(\ref{eq:reward_0}) (RO method) and with the perturbed reward function~(\ref{eq:obj:manipulated}) (CFA method), respectively.
}
\label{fig:waiting}
\end{figure}
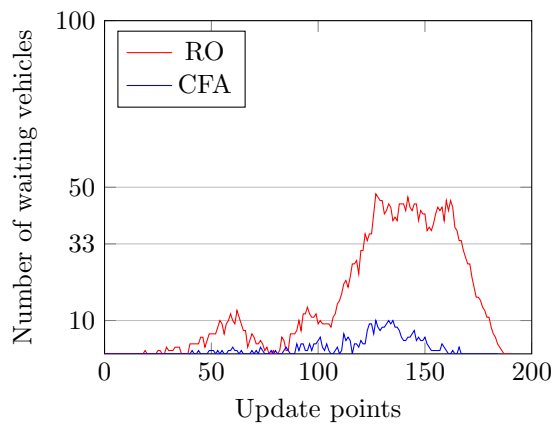
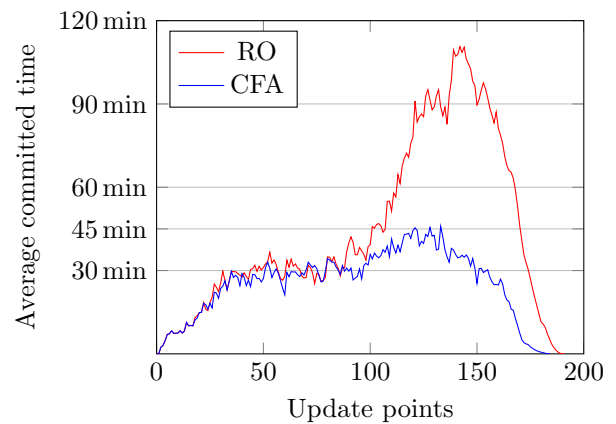

We analyzed the impact of penalizing waiting times on a large problem instance with 4000~orders and 100~vehicles.
First of all, we need a new definition.
The \emph{committed time} of a vehicle at an update time point is the length of the time span that the next decision will no longer affect.
For example, if vehicle~$v$ is at a location at update time~$\tau_k$, then the vehicle is committed until its earliest departure time~$td^{\curr}_{k,v}$, i.e., the committed time is $td^{\curr}_{k,v} - \tau_k$.
If the vehicle is on the way, then it is committed until the end of its service at the destination factory, that is, its committed time is $td^1_{k,v} - \tau_k$.

Consider \Cref{fig:waiting}.
The red plots refer to the case where decision problems were solved with the original reward function~(\ref{eq:reward_0}), and the blue plots refer to the usage of the perturbed reward function~(\ref{eq:obj:manipulated}) with $\lambda_3=0.75\times\lambda_2$.
In \Cref{fig:waiting:vehicles,fig:waiting:time}, we depict for each update point the number of those vehicles that are currently waiting at a location for a free docking port, and the average committed time of the vehicles, respectively.
Both indicators are suitable for measuring the flexibility of routes.

In \Cref{fig:waiting:vehicles}, we can see that without explicitly penalizing the waiting times, there is a long, namely a 45-epoch (7.5~hours) period, when the number of waiting vehicles is between 34 and 48.
That is, on a third of the planning horizon, a third or almost half of the vehicles are waiting.
When the waiting times are penalized in the reward function, the number of waiting vehicles are decreases significantly, namely, never exceeds~10.

The inflexibility of routes in the myopic approach is even more noticeable in the other figure.
According to \Cref{fig:waiting:time}, in some cases the average committed time of the vehicles is more than one and a half hours.
There are vehicles which have to wait more than 6~hours for a free docking port.
Consequently, if those vehicles came to that factory to pickup some orders, then those orders will be delivered at least 2 hours after their due date, inherently.
Explicitly penalizing waiting times also significantly decreases the average committed time, which never exceeds 46~minutes.

\subsubsection{Penalizing the idle vehicles}
We ran our method with objective function $1/n \times \mathbf{f}_1 + \lambda\times \mathbf{f}_2 + \lambda_4\times \mathbf{f}_4$ for all instances.
The multiplier~$\lambda_4$ for the number of idle vehicles was chosen from $\{0, 5, 10\}$.
The average scores for each group, and for all instances are indicated in \Cref{tab:eval:idle}, while detailed results can be found in the supplementary data.

\begin{table}
\caption{Evaluation of penalizing idle vehicles.}
\label{tab:eval:idle}
\begin{threeparttable}
\begin{tabular}{r rrr}
\toprule
& \multicolumn{3}{c}{Multiplier for idle vehicles ($\lambda_4$)}\\
\cmidrule(lr){2-4}
Group & $0$ & $5$ & $10$\\
\midrule
1 &       1\,307.5 &          979.5 &       1\,076.0 \\
2 &      31\,707.6 &      28\,532.8 &      20\,265.7 \\
\addlinespace[1ex]
3 &          690.0 &          887.5 &          887.5 \\
4 &       6\,750.9 &       6\,658.2 &       6\,658.2 \\
\addlinespace[1ex]
5 &      11\,597.1 &       2\,969.8 &       2\,969.8 \\
6 &      49\,011.9 &      21\,919.9 &      21\,919.9 \\
\addlinespace[1ex]
7 &  1\,076\,526.1 &  4\,856\,878.2 &  4\,909\,406.0 \\
8 & 11\,605\,754.6 & 15\,046\,689.2 & 15\,046\,689.2 \\
\cmidrule(lr){1-4}
average: & 1\,597\,918.2 &  2\,495\,689.4 & 2\,501\,234.0 \\
\bottomrule
\end{tabular}
\end{threeparttable}
\end{table}

In contrast to the previous experiments, penalizing idle times has a bad impact on the largest instances (Groups 7-8) but can improve other ones.
The average improvement on instances with 5 vehicles (Groups 1-2) is 10-36\%, and on instances with 50 vehicles (Groups 5-6) is 55-74\%.

\subsubsection{Parameter tuning}
As we can see, the introduced two penalty terms behave differently on the groups.
We tested our algorithm for all pair $(\lambda_3,\lambda_4)$ with respect to the given parameter sets.
We obtained best results (in the sense of total average score) with parameters $\lambda_3 = 0.5 \times \lambda$ and $\lambda_4 = 5$.

Note that since the number of vehicles is part of the input, we could use different parameter settings for the different groups.
However, we do not take advantage of this opportunity in the following, but use the same parameters for all groups.

\subsection{Comparison of existing methods}\label{sec:comp:comp}

\begin{table}
\caption{Comparison of the solution approaches on the benchmark dataset. Best average scores are in bold.}
\label{tab:comp}
\begin{threeparttable}
\begin{tabular}{r r rrr r r}
\toprule
Group & CFA-VNS & Gold & Silver & Bronze & VNSME & MOEA/D-ES\\
\midrule
1 & \textbf{979.5} & 2\,896.4 & 13\,676.2 & 1\,763.8 & 1\,036.8 & 1\,024.2\\
2 & 28\,532.8 & 41\,535.3 & 599\,932.8 & 62\,180.2 & 36\,765.0 & \textbf{14\,182.4}\\
3 & 888.4 & 5\,860.4 & 2\,310.7 & 8\,969.7 & 691.9 & \textbf{686.0}\\
4 & 7\,456.1 & 6\,544.5 & 105\,049.4 & 26\,938.3 & 7\,909.1 & \textbf{5\,489.9}\\
5 & \textbf{3\,159.2} & 10\,459.1 & 17\,284.3 & 94\,794.8 & 10\,541.2 & 9\,673.9\\
6 & \textbf{16\,178.2} & 41\,494.3 & 153\,419.1 & 651\,945.0 & 42\,375.0 & 25\,362.5\\
7 & \textbf{629\,332.6} & 798\,240.7 & 904\,586.3 & 1\,941\,385.3 & 988\,012.5 & 787\,050.0\\
8 & \textbf{4\,470\,954.4} & 11\,359\,466.4 & 18\,678\,529.1 & 15\,122\,816.2 & 10\,337\,500.0 & 10\,260\,000.0\\
\cmidrule(lr){1-7}
average & \textbf{644\,685.2} & 1\,533\,312.1 & 2\,559\,348.5 & 2\,238\,849.2 & 1\,428\,104.0 & 1\,387\,933.6\\
\bottomrule
\end{tabular}
\begin{tablenotes}
\item Gold, Silver, Bronze: Public algorithms of the ICAPS 2021 DPDP Competition
\item VNSME: \citep{cai2022variable}
\item MOEA/D-ES: \citep{cai2022efficient, cai2023survey}
\end{tablenotes}
\end{threeparttable}
\end{table}

In the following, we compare our method to the existing solution approaches.
In \Cref{tab:comp}, we indicate the average group scores of our best method (\textit{CFA-VNS}), the algorithms of the podium teams of the DPDP competition (\textit{Gold}, \textit{Silver}, and \textit{Bronze}), the variable neighborhood search of \citet{cai2022variable} (\textit{VNSME}), and the multi-objective evolutionary method of \citet{cai2022efficient, cai2023survey} (\textit{MOEA/D-ES}).
The first four algorithms (CFA-VNS, Gold, Silver, and Bronze) were tested in the same environment, since the algorithms of the podium teams are publicly available.
Note that the available algorithm of \citet{silversolution} failed on seven instances of Group~2 due to a technical constraint of the simulator, thus the results for those instances are obtained by turning off that constraint. 
Note that methods Gold and VNSME are the same, but the values corresponding to VNSME are obtained from the results in \citep{cai2022variable}.
Method VNSME were executed on a workstation with an Intel Core i5-9500 3.00 GHz CPU with 4 cores, under Ubuntu 18.04 LTS operating system.
According to CPU benchmarks\footnote{\url{https://www.cpubenchmark.net}}, the single thread rating is 2499 for our processor, and 2571 for the other one.
The running environment of method MOEA/D-ES is unknown, the corresponding values are obtained from the results in \citep{cai2023survey}.

Methods CFA-VNS and MOEA/D-ES outperform the other solution approaches on this benchmark dataset.
Comparing these two methods, the multi-objective evolutionary algorithm obtained 26-50\% better average scores on the small instances (Groups 2-4), however, CFA-VNS gave 20-60\% better average scores on the larger instances (Groups 5-8).
The average improvement of CFA-VNS on the full dataset is 54\%.

\begin{table}
\caption{Comparison of the solution approaches on \citeauthor{du2023hierarchical}'s selection.}
\label{tab:comp:du}
\begin{threeparttable}
\begin{tabular}{r rrrr}
\toprule
Instance & CFA-VNS & VNSME\tnote{$\star$} & MOEA/D-ES\tnote{$\star$} & Hierarchical\\
\midrule
38 & \textbf{3\,793.1} & 15\,000.0 & 14\,200.0 & 28\,121.0\\
39 & \textbf{5\,962.1} & 14\,900.0 & 17\,400.0 & 24\,013.0\\
40 & \textbf{7\,049.6} & 10\,200.0 & 13\,100.0 & 23\,338.0\\
\cmidrule(lr){1-5}
41 & \textbf{7\,939.8} & 29\,300.0 & 28\,300.0 & 126\,292.0\\
44 & \textbf{34\,649.0} & 77\,800.0 & 36\,200.0 & 185\,288.0\\
45 & \textbf{10\,244.7} & 36\,300.0 & 22\,800.0 & 185\,909.0\\
\cmidrule(lr){1-5}
49 & 788\,852.0 & 1\,580\,000.0 & \textbf{630\,000.0} & 1\,082\,139.0\\
54 & \textbf{621\,130.3} & 2\,070\,000.0 & 1\,630\,000.0 & 1\,619\,514.0\\
55 & 530\,422.7 & \textbf{323\,000.0} & 1\,000\,000.0 & 678\,125.0\\
\bottomrule
\end{tabular}
\begin{tablenotes}
\item VNSME: \citep{cai2022variable}
\item MOEA/D-ES: \citep{cai2022efficient, cai2023survey}
\item Hierarchical: \citep{du2023hierarchical}
\item[$\star$] values were provided in 2-decimal exponential format (i.e., $E+n$)
\end{tablenotes}
\end{threeparttable}
\end{table}

\citet{du2023hierarchical} provide quantitative results for a selection of 9  instances from the ICAPS 2021 DPDP competition.
Moreover, they compare their results to that of \citet{silversolution}, the silver-winning team of the competition. 
In \Cref{tab:comp:du} we indicate the scores for these selected instances of our best method (\textit{CFA-VNS}), and the hierarchical method of \citet{du2023hierarchical} (\textit{Hierarchical}).
For the sake of completeness, we also indicate the scores for the variable neighborhood search of \citet{cai2022variable} (\textit{VNSME}), and the multi-objective evolutionary method of \citet{cai2022efficient, cai2023survey} (\textit{MOEA/D-ES}). 
Note that the results in the latter papers are provided in a 2-decimal exponential format (e.g., 1.42E+4), thus the values indicated in \Cref{tab:comp:du} for these methods are approximated values.
We can see that CFA-VNS produced better scores on 8 out of 9 instances than VNSME and MOEA/D-ES, respectively.
The difference compared to Hierarchical is even greater, as CFA-VNS gave 22-94\% better results.

Also note that \citet{du2023hierarchical} generated 8 new instances with 50~orders and 5~vehicles, and 8 new instances with 2000~orders and 50~vehicles.
The authors also make comparisons on this new dataset, however, they do not provide any numerical results.
Only bar charts are given, from which it is difficult to read even approximate values.
Since we could not get access to more detailed numerical results so far, we cannot provide  comparison with our solution approach on this dataset.

\section{Conclusions}\label{sec:conc}
In this paper we investigated a dynamic pickup and delivery problem with docking constraints and the LIFO rule.
We proposed a cost function approximation (CFA) method for the problem, where we perturbed the objective function to make solutions flexible for future changes.
The main contribution of our method are the two penalty terms added to the cost function that penalize waiting for service at the locations and also if the vehicles are idle. 
We proposed a variable neighborhood search with three LIFO-specific local search operators to solve problems at the states with the perturbed objective function.
As a case study, we evaluated our solution procedure on the instances of the ICAPS 2021 DPDP competition.
The computational experiments show that our method significantly outperforms the other solution approaches on this dataset.
The average improvement over the state-of-the-art on the full dataset is more than 50\%.
Our method is especially good on the largest (and hardest) instances.

While no exploitable probabilistic information may be available in a single problem instance, statistical data  about the distribution of orders in time and space might be exploited even better than our current method does.
However, we leave it to future work.

\section*{Acknowledgments}
This research has been supported by the TKP2021-NKTA-01 NRDIO grant on "Research on cooperative production and logistics systems to support a competitive and sustainable economy".
Mark\'o Horv\'ath acknowledges the support of the J\'anos Bolyai Research Scholarship.


\bibliographystyle{apalike}
\bibliography{references}


\begin{appendices}

\section{Serving vehicles at factories}\label{sec:apx:serving}
In this section we describe the method of serving the vehicles at a factory in detail.
Recall that each factory has a given number of docking ports for loading and unloading, and the vehicles are served in non-decreasing arrival time order, where ties are resolved randomly.

Given a factory with $c$ ports for loading and unloading the vehicles, and suppose the current time is $t$.
Let $\mathcal{L} = (v_1,v_2,\ldots,v_{k})$ be the \emph{reservation list} consisting of those vehicles which are currently at the factory.
Let $td_i$ be the earliest departure time of vehicle $v_i$.
$\mathcal{L}$ is ordered in non-decreasing departure time order, that is, $td_i \leq td_{i+1}$ for $1\leq i\leq k-1$.
Suppose  vehicle~$v_{k+1}$  arrives at the factory at time $t$.
Finally, let $st_{k+1}$ denote the service time of vehicle~$v_{k+1}$, which is the sum of the dock approaching time, the unloading time and the loading time.

If $k< c$, the vehicle $v_{k+1}$ immediately starts to approach a free port.
The waiting time is zero, and the earliest departure time is $td_{k+1} = t + st_{k+1}$.
Otherwise, if $k \geq c$, then all ports are occupied, and the vehicle must wait until vehicles $v_1,v_2,\ldots,v_{k-c+1}$ finish, that is, until time~$td_{k-c+1}$.
Thus, the waiting time is $td_{k-c+1} - t$, and the earliest departure time is $td_{k+1} = td_{k-c+1} + st_{k+1}$.
In both cases, the vehicle is inserted into the appropriate position of the reservation list.
When a vehicle finishes, it is removed from the list.

\begin{figure}
\begin{tikzpicture}
\tikzstyle{docking}  =[draw,fill=black!0!white];
\tikzstyle{unloading}=[draw,fill=black!25!white];
\tikzstyle{loading}  =[draw,fill=black!50!white];
\tikzstyle{waiting}  =[draw,fill=black];

\def\h{0.5}
\def\va{ 0.00}
\def\vb{-0.75}
\def\vc{-1.50}
\def\vd{-2.25}
\def\t{-2.50}

\node[below] at (0.00,\t) {$t_1$};
\draw[dashed,black!50!white] (0.00,\va) -- (0.00,\t);
\node[below] at (0.75,\t) {$t_2$};
\draw[dashed,black!50!white] (0.75,\vb) -- (0.75,\t);
\node[below] at (1.50,\t) {$t_3$};
\draw[dashed,black!50!white] (1.50,\vc) -- (1.50,\t);
\node[below] at (2.25,\t) {$t_4$};
\draw[dashed,black!50!white] (2.25,\vd) -- (2.25,\t);
\node[below] at (2.75,\t) {$t_5$};
\draw[dashed,black!50!white] (2.75,\vb) -- (2.75,\t);
\node[below] at (3.25,\t) {$t_6$};
\draw[dashed,black!50!white] (3.25,\va) -- (3.25,\t);
\node[below] at (5.00,\t) {$t_7$};
\draw[dashed,black!50!white] (5.00,\vc) -- (5.00,\t);
\node[below] at (6.00,\t) {$t_8$};
\draw[dashed,black!50!white] (6.00,\vd) -- (6.00,\t);

\draw[thick,->,-latex] (-0.25,\t) -- (6.50,\t);

\node[left] at (-0.25,\va+0.5*\h) {vehicle 1};
\node[left] at (-0.25,\vb+0.5*\h) {vehicle 2};
\node[left] at (-0.25,\vc+0.5*\h) {vehicle 3};
\node[left] at (-0.25,\vd+0.5*\h) {vehicle 4};

\draw[docking]   (0.00,\va) rectangle (1.00,\va+\h);
\draw[unloading] (1.00,\va) rectangle (2.50,\va+\h);
\draw[loading]   (2.50,\va) rectangle (3.25,\va+\h);

\draw[docking]   (0.75,\vb) rectangle (1.75,\vb+\h);
\draw[unloading] (1.75,\vb) rectangle (2.25,\vb+\h);
\draw[loading]   (2.25,\vb) rectangle (2.75,\vb+\h);

\draw[waiting]   (1.50,\vc + 0.25*\h) rectangle (2.75,\vc + 0.75*\h);
\draw[docking]   (2.75,\vc) rectangle (3.75,\vc+\h);
\draw[unloading] (3.75,\vc) rectangle (4.50,\vc+\h);
\draw[loading]   (4.50,\vc) rectangle (5.00,\vc+\h);

\draw[waiting]   (2.25,\vd + 0.25*\h) rectangle (3.25,\vd + 0.75*\h);
\draw[docking]   (3.25,\vd) rectangle (4.25,\vd+\h);
\draw[unloading] (4.25,\vd) rectangle (4.75,\vd+\h);
\draw[loading]   (4.75,\vd) rectangle (6.00,\vd+\h);

\node[right] at (3.50,\va+0.5*\h) {\textit{port~1}};
\node[right] at (3.50,\vb+0.5*\h) {\textit{port~2}};
\node[right] at (5.25,\vc+0.5*\h) {\textit{port~2}};
\node[right] at (6.25,\vd+0.5*\h) {\textit{port~1}};
\end{tikzpicture}
\caption{
  Example for vehicles arriving at a factory with two docking ports.
  White rectangles represent dock approaching, gray rectangles represent unloading and loading orders, black rectangles represent waiting for ports to become free.
}
\label{fig:ex:ports}
\end{figure}
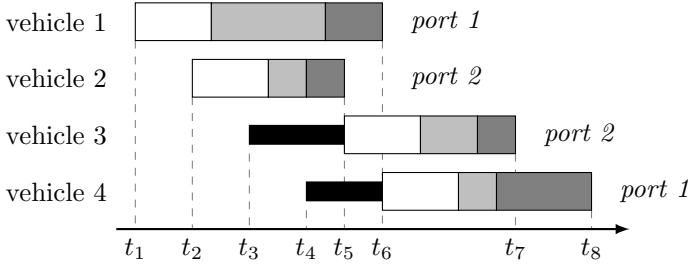

\paragraph{Example}
In \Cref{fig:ex:ports} we depict a situation where four vehicles arrive at a factory with two docking ports.
Vehicle~$v_1$ arrives at the factory at time~$t_1$ and occupies a free docking port ($td_1 = t_6$, $\mathcal{L} = (v_1)$).
Vehicle~$v_2$ arrives at time~$t_2$ and occupies the other free port ($td_2 = t_5$, $\mathcal{L} = (v_2, v_1)$).
Vehicle~$v_3$ arrives at time~$t_3$, however, since both ports are in use (that is, currently vehicle~$v_1$ is under unloading at the first port, and vehicle~$v_2$ approaches the other one), it has to wait until a port is freed at time~$t_5$ ($td_3 = t_7$, $\mathcal{L} = (v_2,v_1,v_3)$).
Vehicle~$v_4$ arrives at time~$t_4$, however, both ports are in use, moreover, vehicle~$v_3$ is already allocated to the port becoming free at time~$t_5$, thus, it has to wait until time~$t_6$, when loading is finished at the first port ($td_4 = t_8$, $\mathcal{L} = (v_2,v_1,v_3,v_4)$).

\section{Feasibility of route plans in states}\label{sec:apx:feasibility}
In this section we characterize feasible  route plans in some state $s_k$.

\paragraph{Orders}
For a vehicle~$v$, let $\mathcal{O}_{k,v}$ be the set of those orders that are carried by the vehicle or picked up in its route plan $\theta_{k,v}$.
That is, $\mathcal{O}_{k,v} = \mathcal{C}_{k,v} \cup \bigcup_{j=1}^{\ell_{k,v}} \Pk^j_{k,v}$.
These, and only these orders must be delivered by the vehicle in the route plan, that is,
\[
\bigcup_{j=1}^{\ell_{k,v}} \Dl^j_{k,v} = \mathcal{O}_{k,v}.
\]
Moreover, the sets $\Pk^j_{k,v}$,  $\Pk^{j'}_{k,v}$ must be disjoint for $j\neq j'$, and $\mathcal{C}_{k,v}$ must be disjoint from $\bigcup_{j=1}^{\ell_{k,v}} \Pk^j_{k,v}$. Further on, 
$\mathcal{C}_{k,v} \cup \bigcup_{j=1}^{\ell_{k,v}} \Pk^j_{k,v}$ must be disjoint from $\mathcal{C}_{k,w} \cup \bigcup_{j=1}^{\ell_{k,w}} \Pk^j_{k,w}$ for distinct vehicles $v \neq w$.
Clearly, orders must be picked up before their delivery, that is, if $o_i \in \Pk^{j_1}_{k,v}$ for some $j_1$, then $o_i \in \Dl^{j_2}_{k,v}$ for some $j_1 < j_2$.
Finally, $\bigcup_{j=1}^{\ell_{k,v}} \Pk^j_{k,v}$ must be a subset of $\tilde{\omega}_k$ for each $v \in \mathcal{V}$.

\paragraph{LIFO constraint}
Let $\mathcal{L}_{k,v}$ be the concatenation of lists $\mathcal{C}_{k,v}$, $\mathcal{D}^1_{k,v}$, $\mathcal{P}^1_{k,v}$, $\ldots$, $\mathcal{D}^{\ell_{k,v}}_{k,v}$, $\mathcal{P}^{\ell_{k,v}}_{k,v}$.
Let $\operatorname{pos}(o_i^+)$ and $\operatorname{pos}(o_i^-)$ denote the position of the first and the second (i.e., last) occurrence of order $o_i \in \mathcal{O}_{k,v}$ in $\mathcal{L}_{k,v}$, respectively.
Then, route plan~$\theta_{k,v}$ satisfies the LIFO constraint, if
\[
  \operatorname{pos}(o_i^+) < \operatorname{pos}(o_j^+) \Rightarrow \operatorname{pos}(o_i^-) \leq \operatorname{pos}(o_j^+) \vee \operatorname{pos}(o_j^-) \leq \operatorname{pos}(o_i^-)
\]
holds for all $o_i, o_j \in \mathcal{O}_{k,v}$.

\paragraph{Capacity constraint}
The route plan~$\theta_{k,v}$ satisfy the capacity constraint, if the total quantity of the loaded orders never exceeds the vehicle's capacity.
That is,
\[
  \sum_{o_i\in\mathcal{C}_{k,v}} q_i + \sum_{j=1}^{\ell} \left( \sum_{o_i\in \mathcal{P}_{k,v}^j} q_i - \sum_{o_i\in \mathcal{D}_{k,v}^j} q_i \right) \leq Q
\]
holds for all $\ell=1,\ldots,\ell_{k,v}$.

\paragraph{Fundamental routing constraints}
The travel time between factories are fixed:
\[
ta_v^j-td_v^{j-1} = \operatorname{travel}(f_v^{j-1},f_v^j)\quad \text{for all}\ j=1,\ldots,\ell_v.
\]
Let $\eta_v^j$ be the waiting time of the vehicle at the $j$th visited factory.
\[
  ta_v^j + \eta_v^j + h^{\operatorname{docking}} + \sum_{o_k \in \mathcal{D}_v^j} h^d_k + \sum_{o_k \in \mathcal{P}_v^j} h^p_k \leq td_v^j
\]

\section{Transition}\label{sec:apx:transition}
In the following, we formally describe the transition from state~$s_k$ to state~$s_{k+1}$ according to action~$x_k$, postponed from \Cref{sec:prob:transition}.
The various cases are summarized in \Cref{tab:transition}, and explained in the following.
Recall that $\theta^i_{k,v}$ refers to the $i$th visit of vehicle~$v$ in its route plan belongs to state~$s_k$, and $\theta^{x,i}_{k,v}$ refers to the $i$th visit of the route plan belongs to action~$x_k$.

\begin{table}
\caption{Transition from state~$s_k$ into state~$s_{k+1}$ according to action~$x_k$.}
\label{tab:transition}
\begin{threeparttable}
\begin{tabular}{cccccccc}
\toprule
& \multicolumn{2}{c}{state~$s_k$} & \multicolumn{2}{c}{action~$x_k$} & \multicolumn{3}{c}{state~$s_{k+1}$}\\
\cmidrule(lr){2-3} \cmidrule(lr){4-5} \cmidrule(lr){6-8}
Case & $\phi^{\curr}_{k,v}$ & $\theta^1_{k,v}$ & $\theta^{x,1}_{k,v}$ & $\theta^{x,2}_{k,v}$ & $\tau_{k+1}$ & $\phi^{\curr}_{k+1,v}$ & $\theta^1_{k+1,v}$ \\
\midrule
1 & yes & $\star$ & $\star$ & $\star$ & $\tau_{k+1} < td^{\curr}_{k,v}$ & $\phi^{\curr}_{k,v}$& $\theta^1_{k,v}$ \\
2a & yes & $\star$ & no & no & $td^{\curr}_{k,v} \leq \tau_{k+1}$ & $(f^{\curr}_{k,v},\tau_{k+1})$ & no \\
2b & yes & $\star$ & yes & $\star$ & $td^{\curr}_{k,v} \leq \tau_{k+1}$ & no & $\theta^{x,1}_{k,v}$ \\
3 & no & yes & yes & $\star$ & $\tau_{k+1} < ta^1_{k,v}$ & no & $\theta^1_{k,v}$ \\
4 & no & yes & yes & $\star$ & $ta^1_{k,v} \leq \tau_{k+1}$ & $(f^1_{k,v},td^1_{k,v})$ & $\theta^{x,2}_{k,v}$ \\
\bottomrule
\end{tabular}
\begin{tablenotes}
\item[$\star$] could be 'yes' or 'no' (i.e., 'given' or 'not given').
Recall that if $\theta^1_{k,v}$ is given, so is $\theta^{x,1}_{k,v}$ with $f^1_{k,v} = f^{x,1}_{k,v}$ and $ta^1_{k,v} = ta^{x,1}_{k,v}$.
\end{tablenotes}
\end{threeparttable}
\end{table}

\paragraph{Case 1. The vehicle is not finished at its current factory}
If the vehicle had a current factory and the associated earliest departure time is later then the current update point, then the current factory remains the same as in the previous state, and the destination is the one specified in the previous action, if any.

\paragraph{Case 2. The vehicle is finished at its current factory}
Assume that the vehicle had a current factory and the associated earliest departure time has passed.
If the vehicle was not assigned a destination factory in the previous action (Case 2a), then the vehicle remains parked at this factory, but will be immediately available.
Otherwise (Case 2b), the vehicle is already on the way to its destination specified in the last action.

\paragraph{Case 3. The vehicle is not reached its destination factory}
If a vehicle was on the way to its destination factory, and did not reach it in the last epoch, then the destination remains the same (however, the list of orders to pickup may be changed in the action).

\paragraph{Case 4. The vehicle has reached its destination factory}
Assume that the vehicle was on the way to its destination factory, and reached that in the last epoch.
Then, that former destination becomes the current factory, and earliest departure time is also calculated.
If the vehicle was not assigned further factories to visit in the previous action, then the vehicle is not assigned a destination factory in the current state.
Otherwise, the first factory to visit becomes the next destination factory.

\paragraph{Carrying orders}
In all cases $\Pk^0_{k+1,v} = \Pk^0_{k,v}$, except in Case~4, when $\Pk^0_{k+1,v} = (\Pk^0_{k,v} \setminus \mathcal{D}^1_{k,v}) \cup \mathcal{P}^1_{k,v}$, since the vehicle has reached its destination, although loading and unloading may be in progress at time point $\tau_{k+1}$.

\section{Feasibility of solutions}\label{sec:apx:feas_solution}
In the following, we formally define the feasibility of the route plans in a solution.

\paragraph{Orders}
First of all, all orders must be dispatched in a solution.
Let $\mathcal{O}_v$ be the set of those orders which belong to (i.e., picked up and delivered by) vehicle~$v$.
Then, each order belongs to exactly one vehicle:
\[
  \bigcup_{v\in \mathcal{V}} \mathcal{O}_v = \mathcal{O} 
  ,
\]
with $\mathcal{O}_v \cap \mathcal{O}_w = \emptyset$ for all $v\neq w$.
Also, each order is picked up and delivered only once:
\[
  \mathcal{O}_v = \bigcup_{j=1}^{\ell_v} \mathcal{D}_v^j = \bigcup_{j=1}^{\ell_v} \mathcal{P}_v^j,
\]
with $\mathcal{D}_v^i\cap \mathcal{D}_v^j = \emptyset$ and $\mathcal{P}_v^i\cap \mathcal{P}_v^j = \emptyset$ for all $i\neq j$.

Each order must be picked up before its delivery, that is, if $o_i \in \Pk^{j_1}_{v}$ for some $j_1$, then $o_i \in \Dl^{j_2}_{v}$ for some $j_1 < j_2$.

\paragraph{LIFO constraint}
Let $\mathcal{L}_{v}$ be the concatenation of lists $\mathcal{D}^1_{v}$, $\mathcal{P}^1_{v}$, $\ldots$, $\mathcal{D}^{\ell_{v}}_{v}$, $\mathcal{P}^{\ell_{v}}_{v}$.
Let $\operatorname{pos}(o_i^+)$ and $\operatorname{pos}(o_i^-)$ denote the position of the first and the second (i.e., last) occurrence of order $o_i \in \mathcal{O}_{v}$ in $\mathcal{L}_{v}$, respectively.
Then, route plan~$\theta_{v}$ satisfies the LIFO constraint, if
\[
  \operatorname{pos}(o_i^+) < \operatorname{pos}(o_j^+) \Rightarrow \operatorname{pos}(o_i^-) \leq \operatorname{pos}(o_j^+) \vee \operatorname{pos}(o_j^-) \leq \operatorname{pos}(o_i^-)
\]
holds for all $o_i, o_j \in \mathcal{O}_{v}$.

\paragraph{Capacity constraint}
The route plan~$\theta_v$ satisfy the capacity constraint, if the total quantity of the orders carried by the vehicle never exceeds its capacity.
That is,
\[
  \sum_{j=1}^\ell \left( \sum_{o_i\in \mathcal{P}_v^j} q_i - \sum_{o_i\in \mathcal{D}_v^j} q_i \right) \leq Q
\]
holds for all $\ell=1,\ldots,\ell_v$.

\paragraph{Fundamental routing constraints}
The travel time between factories are fixed:
\[
ta_v^j-td_v^{j-1} = \operatorname{travel}(f_v^{j-1},f_v^j)\quad \text{for all}\ j=1,\ldots,\ell_v.
\]
Let $\eta_v^j$ be the waiting time of the vehicle at the $j$th visited factory.
\[
  ta_v^j + \eta_v^j + h^{\operatorname{docking}} + \sum_{o_k \in \mathcal{D}_v^j} h^d_k + \sum_{o_k \in \mathcal{P}_v^j} h^p_k \leq td_v^j
\]

\section{Evaluation procedure}
\label{sec:apx:eval}
The purpose of the evaluation procedure is to compute the value of the cost functions (\ref{eq:obj}) and (\ref{eq:obj:manipulated}).
To this end, the timing information for the routes defined in \Cref{sec:method} has to be computed.
We apply a basic event-based simulation procedure.
The procedure maintains a priority queue (\emph{event queue}) in order to process arrival and departure events.
In addition, each factory is associated with a \emph{reservation list} containing vehicles that are currently at the factory, see Appendix~\ref{sec:apx:serving}.

\paragraph{Step 1 (Initialization)}
First, consider the vehicles with current factory.
If the earliest departure time of a vehicle is in the future (i.e., the vehicle is not finished yet at this factory), we insert the vehicle into the sorted reservation list of the factory.
We also create a \emph{departure event} associated with the corresponding departure time.
Then, for each vehicle without current factory, we create an \emph{arrival event} associated with the corresponding arrival time.

\paragraph{Step 2 (Processing events)}
If the event queue is empty, we are done, and thus stop by returning the calculated objective function value.
Otherwise, we take out from the queue an event with the earliest associated time.
Denote the corresponding factory with~$f_e$, and the associated time with~$t_e$.
If this event is a departure event, we remove the corresponding vehicle from the reservation list of the corresponding factory.
Also, if the vehicle has a next factory to visit, say $f_n$, we create and store an arrival event associated with the arrival time $t_e + \operatorname{travel}(f_e,f_n)$.
Otherwise, if the event taken out is an arrival event, then we calculate and store its contribution to the cost function value, as defined in \Cref{sec:prob:dynamic:obj}.
We also calculate the earliest departure time of the vehicle (see Appendix~\ref{sec:apx:serving}), and create a departure event associated with this time.
We repeat this step.
\end{appendices}
\end{document}